\documentstyle[amscd,amssymb,verbatim,diagrams,12pt]{amsart}
\newcommand{\FM}{\operatorname{FM}}
\newcommand{\bH}{{\bf H}}
\newcommand{\BU}{{\bf U}}
\newcommand{\BT}{{\bf T}}
\newcommand{\BZ}{{\bf Z}}
\newcommand{\Bf}{{\bf f}}
\newcommand{\ZZ}{{\cal Z}}
\newcommand{\Id}{\operatorname{Id}}

\newcommand{\HP}{{\frak H}}
\newcommand{\ch}{\operatorname{ch}}

\renewcommand{\mod}{\operatorname{mod}}

\newcommand{\Tr}{\operatorname{Tr}}

\newcommand{\OO}{{\cal O}}

\newcommand{\NN}{{\cal N}}

\renewcommand{\Re}{\operatorname{Re}}
\newcommand{\SL}{\operatorname{SL}}
\newcommand{\G}{{\Bbb G}}

\newcommand{\hra}{\hookrightarrow}
\newcommand{\lan}{\langle}
\newcommand{\ran}{\rangle}
\newcommand{\Coh}{\operatorname{Coh}}

\newcommand{\CC}{{\cal C}}

\renewcommand{\Sp}{\operatorname{Sp}}
\newcommand{\Mat}{\operatorname{Mat}}

\renewcommand{\P}{{\Bbb P}}

\newcommand{\si}{\sigma}

\newcommand{\Pic}{\operatorname{Pic}}

\newcommand{\ga}{\gamma}
\newcommand{\de}{\delta}
\newcommand{\eps}{\epsilon}

\renewcommand{\ker}{\operatorname{ker}}

\numberwithin{equation}{subsection}

\newcommand{\GL}{\operatorname{GL}}

\newcommand{\Stab}{\operatorname{Stab}}

\newtheorem{thm}{Theorem}[subsection]
\newtheorem{prop}[thm]{Proposition}
\newtheorem{lem}[thm]{Lemma}
\newtheorem{cor}[thm]{Corollary}
{  \theoremstyle{definition}
\newtheorem{defi}[thm]{Definition}
\newtheorem{ex}[thm]{Example}
\newtheorem{exs}[thm]{Examples}
\newtheorem{rem}[thm]{Remark}
\newtheorem{rems}[thm]{Remarks}
}

\newcommand{\Pf}{\noindent {\it Proof}}
\newcommand{\id}{\operatorname{id}}

\newcommand{\ov}{\overline}

\renewcommand{\Im}{\operatorname{Im}}

\newcommand{\rk}{\operatorname{rk}}
\newcommand{\ra}{\rightarrow}
\renewcommand{\AA}{{\cal A}}
\newcommand{\Ab}{{\cal A}b}
\newcommand{\FF}{{\cal F}}

\newcommand{\TT}{{\cal T}}

\newcommand{\PP}{{\cal P}}

\newcommand{\LL}{{\cal L}}

\newcommand{\Om}{\Omega}

\newcommand{\St}{\operatorname{St}}
\newcommand{\Hom}{\operatorname{Hom}}

\newcommand{\Ext}{\operatorname{Ext}}
\newcommand{\End}{\operatorname{End}}

\renewcommand{\a}{\alpha}
\renewcommand{\b}{\beta}
\newcommand{\om}{\omega}
\newcommand{\De}{\Delta}
\newcommand{\la}{\lambda}

\newcommand{\C}{{\Bbb C}}
\newcommand{\N}{{\Bbb N}}
\newcommand{\R}{{\Bbb R}}
\newcommand{\Z}{{\Bbb Z}}
\newcommand{\Q}{{\Bbb Q}}
\newcommand{\La}{\Lambda}
\newcommand{\Ga}{\Gamma}

\newcommand{\wt}{\widetilde}
\newcommand{\ot}{\otimes}
\newcommand{\sign}{\operatorname{sign}}

\newcommand{\sub}{\subset}
\newcommand{\ed}{\qed\vspace{3mm}}

\newcommand{\SH}{\operatorname{SH}}
\newcommand{\LG}{{\bf LG}}

\newcommand{\Fun}{\operatorname{Fun}}
\newcommand{\NS}{\operatorname{NS}}
\newcommand{\Arg}{\operatorname{Arg}}
\newcommand{\Pol}{\operatorname{Pol}}
\newcommand{\Spin}{\operatorname{Spin}}
\newcommand{\spin}{\operatorname{spin}}
\newcommand{\linspan}{\operatorname{span}}
\newcommand{\PGL}{\operatorname{PGL}}
\newcommand{\Cl}{\operatorname{Cl}}
\newcommand{\BN}{{\bf N}}
\newcommand{\BP}{{\bf P}}

\newcommand{\Nm}{\operatorname{Nm}}
\newcommand{\cha}{\operatorname{char}}

\newcommand{\bv}{{\bf v}}

\title[Phases of Lagrangian-invariant objects on an abelian variety]{Phases of Lagrangian-invariant objects in the derived category of an abelian variety}
\author{Alexander Polishchuk}
\address{Department of Mathematics, University of Oregon, Eugene, OR 97405}
\email{apolish@@uoregon.edu}
\thanks{Supported in part by NSF grant}

\begin{document}
\begin{abstract} We continue the study of
Lagrangian-invariant objects (LI-objects for short) in the derived
category $D^b(A)$ of coherent sheaves on an abelian variety, 
initiated in \cite{P-LIF}. For every element of the complexified ample cone $D_A$
we construct a natural phase function on the set of LI-objects, which in 
the case $\dim A=2$ gives the phases with respect to 
the corresponding Bridgeland stability (see \cite{Bridge-K3}). 
The construction is based on the relation between endofunctors of $D^b(A)$
and a certain natural central extension of groups, associated with $D_A$ viewed as a hermitian symmetric space. In the case when $A$ is a power of an elliptic curve, we show
that our phase function has a natural interpretation in terms of the Fukaya category of the mirror dual
abelian variety. As a byproduct of our study of LI-objects we show that the Bridgeland's component of
the stability space of an abelian surface contains all full stabilities.
\end{abstract}
\maketitle

\bigskip

\centerline{\sc Introduction}

The notion of stability condition on triangulated categories, introduced by Bridgeland in \cite{Bridge-stab}, axiomatizes the notion of stability of branes coming from the study of deformations of
superconformal field theories (see \cite{Douglas}). The hope is that the space of stability conditions
on a Calabi-Yau threefold is related to the moduli spaces of complex structures on a mirror
dual manifold. At present we have examples of Bridgeland stabilities on $D^b(X)$ for any surface $X$,
however, the problem of constructing such examples for a Calabi-Yau threefold is still open
(see \cite{BMT} for a proposal of such a construction). 

The goal of this paper is to test the existence of a stability condition on $D^b(A)$ for any abelian variety $A$ by looking at certain special objects in $D^b(A)$. More precisely,
for an element $\om=i\a+\b\in D_A\sub\NS(A)\ot\C$ in the complexified ample cone (defined by the condition that $\a$ is ample) one expects to have a stability
condition on $D^b(A)$ with the central charge
$$Z(F)=\int_A\exp(-\om)\cdot\ch(F),$$
where $F\in D^b(A)$. The starting point of this work is the observation that there are certain objects in $D^b(A)$ that are automatically semistable with respect to any nice stability condition (see
Prop.\ \ref{LI-semistable-prop}). Namely, these are {\it Lagrangian-invariant
objects} ({\it LI-objects} for short) defined in \cite{P-LIF} (see also Def.\ \ref{LI-defi}).
The simplest examples are the structure sheaves of points $\OO_x$.
To get other examples one can consider images of $\OO_x$ under autoequivalences of $D^b(A)$
but in general these do not exhaust all LI-objects (see Remark \ref{finite-rem} and
Prop.\ \ref{class-prop}).
Thus, a stability condition should give a {\it phase} for any LI-object $F$, i.e., a lifting
of $\Arg Z(F)\in \R/2\pi\Z$ to $\R$. Furthermore, a nonzero morphism $F_1\to F_2$ can exist
only if the phase of $F_1$ does not exceed the phase of $F_2$. The main result of this paper
is the construction of such a phase function associated with each $\om\in D_A$. We also
verify some properties of this function that one expects from the theory of stability conditions
(see Thm.\ \ref{phase-thm}).

The major role in our construction is played by the action of a certain group on the set 
$\ov{\SH}^{LI}/\N^*$ of
classes of LI-objects modulo certain simple equivalence relation (we allow to apply translations and
tensoring with a line bundle in $\Pic^0(A)$ and with a vector space). This group, which we
denote $\wt{\BU(\Q)}$ is a central extension by $\Z$ of the group of $\Q$-points of
an algebraic group $\BU=\BU_{X_A}$ defined as automorphisms of the abelian variety $X_A=A\times\hat{A}$, compatible with the skew-symmetric autoduality of $X_A$.
The preimage of the subgroup of $\Z$-points 
in $\wt{\BU(\Q)}$ is closely related to the group of autoequivalences of $D^b(A)$
(see \cite{Mukai-spin,P-sym,Orlov-ab}). The main idea that brings the Siegel domain $D_A$
into picture is that the above central extension has a natural interpretation in terms of the
action of $\BU(\Q)$ on $D_A$ (see Theorem \ref{extension-thm}). This allows us to parametrize
the set $\ov{\SH}^{LI}/\N^*$ of classes of LI-objects by points of a natural $\Z$-covering of
the set of $\Q$-points of a certain homogeneous algebraic variety $\LG=\LG_A$ for the group $\BU$
(the points $\LG(\Q)$ are in bijection with {\it Lagrangian} abelian subvarieties in $A\times\hat{A}$),
and the phase function appears naturally in this context.

If $\dim A=2$ then the stability condition corresponding to $\om$ was
constructed by Bridgeland in \cite{Bridge-K3}, and we check that our phases for LI-objects 
match the ones coming from this stability condition (see Section \ref{surface-sec}).

In the case when $A=E^n$, where $E$ is an elliptic curve without complex multiplication,
we give a mirror-symmetric interpretation of our picture in terms of Fukaya category of the mirror dual
abelian variety (following the recipe of \cite{GLO}). We show that the central charge on LI-objects in $D^b(A)$ defined using $\om\in D_A$ matches with the integral of the holomorphic volume form
over the corresponding Lagrangian tori, and hence, that LI-objects in $D^b(A)$ give rise to graded Lagrangians on the mirror dual side (see Section \ref{mirror-sec}).

We also observe that the set $\wt{\LG(\Q)}$ parametrizing classes of LI-objects also parametrizes
certain natural collection of $t$-structures on $D^b(A)$, generalizing the ones obtained from the
standard $t$-structure by applying autoequivalences (we call them {\it quasi-standard}). 
We conjecture that there is also a natural $t$-structure associated with every point of
$\wt{\LG(\R)}$ whose heart is equivalent to the category of holomorphic bundles on the 
corresponding noncommutative torus (see \cite{PS, P-nc-tori, Block}). 

Another by-product of our study is a refinement of the results of \cite{Mukai-spin,GLO} on
the action of autoequivalences of $D^b(A)$ on numerical classes of objects.
Namely, we construct a natural double covering 
$\Spin\to \BU$ of algebraic groups over $\Q$ and an algebraic representation of
$\Spin$ on the vector space associated with the numerical Grothendieck group of $A$,
such that the action of elements projecting to $\BU(\Q)$ is induced by endofunctors of $D^b(A)$
(see Thm.\ \ref{alg-spin-action-thm}).

The paper is organized as follows. Section \ref{prelim-sec} contains some auxiliary results
not involving derived categories. In particular, we give an interpretation of the index of
a nondegenerate line bundle on an abelian variety in terms of the function $\Arg\chi$
on the complexified ample cone (see Theorem \ref{index-thm}). We also prove some
useful results about the group $\BU$ and the variety of Lagrangian subvarieties $\LG$
in $A\times\hat{A}$. In Section \ref{LI-sec} we study the central extension $\wt{\BU(\Q)}\to\BU(\Q)$
coming from a natural $1$-cocycle with values in $\OO^*(D_A)$ and its action on
LI-objects in $D^b(A)$ and their numerical classes. In Section \ref{stab-sec}
we parametrize LI-objects (up to certain equivalence) by points of a natural $\Z$-covering 
$\wt{\LG(\Q)}\to\LG(\Q)$ equipped with an action
of $\wt{\BU(\Q)}$, and construct a family of phase functions on $\wt{\LG(\Q)}$
parametrized by $D_A\times\C$, equivariantly with respect to $\wt{\BU(\Q)}$.
We also study the connection with Bridgeland stability conditions
on abelian surfaces (see Thm.\ \ref{surface-thm}) and with mirror symmetry 
(see Sec.\ \ref{mirror-sec}). Also, as a byproduct of our study of LI-objects we show that any full stability
on an abelian surface $A$ belongs to the Bridgeland's component of the stability space $\Stab(A)$
(see Prop.\ \ref{ab-sur-stab-prop}).
In Section \ref{t-str-FM-sec} we construct a family of $t$-structures on $D^b(A)$ parametrized
by $\wt{\LG(\Q)}$ and study a relation between $\LG(\Q)/\BU(\Q)$ and the Fourier-Mukai partners of $A$ (see Sec.\ \ref{FM-sec}).

\medskip

\noindent
{\it Notations and conventions}. We work over a fixed algebraically closed field $k$.
We say that an object $F$ of a $k$-linear category $\CC$ is {\it endosimple} if
$\Hom_\CC(F,F)=k$.
For a scheme $X$ we denote by $D^b(X)$ the bounded derived category of coherent sheaves on $X$.
We say that an object $F\in D^b(X)$ is {\it cohomologically pure} if there exists a coherent sheaf $H$
such that $F\simeq H[n]$ for $n\in\Z$. 
We denote by $\Ab_\Q$ the category of abelian varieties up to an isogeny (i.e., the localization
of the category of abelian varieties over $k$ with respect to the class of isogenies).
When we want to consider the $F$-vector space associated with a $\Z$-lattice $M$, 
where $F=\Q,\R,$ or $\C$, as an algebraic variety over $F$,
we denote it by $M_F$.

\medskip

\noindent
{\it Acknowledgment}. I am grateful to Tom Bridgeland for helpful discussions
and to Maxim Kontsevich for a discussion of the picture in Section 
\ref{mirror-sec} involving mirror symmetry.


\section{Preliminaries}\label{prelim-sec}

Throughout this paper $A$ denotes an abelian variety over $k$.


\subsection{Degree, trace and Euler bilinear form}

Recall that for $f\in\End(A)$ one has
$$\deg(f)=\det T_l(f)$$
where $T_l(f)$ is the representation of $f$ on the Tate module $T_l(A)$ for $l\neq \cha(k)$
(see \cite[Ch.\ 19, Thm.\ 4]{Mum-ab}).
Thus, extending $\deg$ to a polynomial function
$$\deg: \End(A)\ot\Q\to\Q,$$
homogeneous of degree $2g$, we have 
$$\deg(1+tf)=1+t\cdot\Tr(f)+O(t^2),$$
where $\Tr(f)$ is given by the trace of the action of $f$ on $T_l(A)\ot_{\Z_l}\Q_l$.
Furthermore, $f\mapsto \Tr(f\cdot f')$ is a positive definite quadratic form 
on $\End(A)\ot\Q$, where $f\mapsto f'$ is the Rosati involution associated with a
polarization of $A$ (see \cite[Ch.\ 21, Thm.\ 1]{Mum-ab}).

Let us fix a polarization on $A$ and denote by $\End(A)^+\ot\Q\sub\End(A)\ot\Q$ the subspace of elements invariant with respect to the corresponding Rosati involution.
Note that the quadratic form $\Tr(f^2)$ on $\End(A)^+\ot\Q$ is positive-definite.

\begin{prop}\label{deg-inj-prop}
An element $f\in\End(A)\ot\C$ is determined by the polynomial function
$$\End(A)\ot\C\to\C: x\mapsto \deg(f-x).$$
Furthermore, if $f$ is invariant with respect to the Rosati involution then
it is determined by the restriction of the above function to $\End(A)^+\ot\C$.
\end{prop}

\Pf . We have to check that if $\deg(f_1-x)=\deg(f_2-x)$ for all $x\in\End(A)$ then
$f_1=f_2$. Adding to $f_1$ and $f_2$ the same element of $\End(A)\ot\C$ we can
assume that $f_1$ and $f_2$ are invertible in $\End(A)\ot\C$.
Observe also that $\deg(f_1)=\deg(f_2)$ (this follows by substituting $x=0$).
Thus, we obtain
$$\deg(1-xf_1^{-1})=\deg(f_1-x)\deg(f_1)^{-1}=\deg(f_2-x)\deg(f_2)^{-1}=\deg(1-xf_2^{-1}).$$
Considering the linear terms in $x$ we derive
$$\Tr(x f_1^{-1})=\Tr(x f_2^{-1}).$$
The nondegeneracy of the form $\Tr(fg)$ implies $f_1^{-1}=f_2^{-1}$.

To prove the second statement, we repeat the above argument letting $x$ vary only in
$\End(A)^+\ot\C$. 
\ed

We always use the standard identification 
$$\NS(A)\ot\Q\rTo{\sim}\Hom(A,\hat{A})^+\ot\Q: L\mapsto \phi_L,$$
where $\Hom(A,\hat{A})^+\ot\Q\sub \Hom(A,\hat{A})\ot\Q$ consists
of self-dual homomorphisms. 
The Euler characteristic defines a polynomial function $\chi:\NS(A)\ot\C\to\C$ of degree $g=\dim A$,
which we also view as a function on $\Hom(A,\hat{A})^+\ot\C$.
One has $\chi^2=\deg$ (see \cite[ch.\ 16]{Mum-ab}).

Recall that the Grothendieck group $K_0(A)$ carries the Euler bilinear form 
$$\chi([E],[F]):=\sum_i (-1)^i\dim\Hom^i(E,F),$$
where $E,F\in D^b(A)$. We denote by $\NN(A)$ the numerical Grothendieck group, i.e.,
the quotient of $K_0(A)$ by the kernel of this form. 
$\NN(A)$ is a free abelian group of finite rank (see \cite[Ex.\ 19.1.4]{Fulton}). 
Associating with a line bundle $L$ its class $[L]$ in
$\NN(A)$ defines a polynomial map between free abelian groups of finite rank
$$\ell:\NS(A)\to\NN(A).$$
Therefore, we have the induced polynomial morphism between the corresponding
$\Q$-vector spaces
\begin{equation}\label{ell-map-eq}
\ell:\NS(A)_\Q\to\NN(A)_\Q.
\end{equation}

\begin{cor}\label{chi-cor} An element $\phi\in\NS(A)\ot\C$ is determined by the corresponding
polynomial function 
$$\NS(A)\to\C: x\mapsto \chi(\ell(\phi),\ell(x)).$$
\end{cor}

\Pf . Since $\NS(A)$ is Zariski dense in $\NS(A)_\C$, it is enough to prove the similar statement
with the polynomial function $\chi(\ell(\phi),\ell(\cdot))$ on $\NS(A)\ot\C$.
Note that 
$$\chi(\ell(\phi),\ell(x))^2=\chi(\ell(x-\phi))^2=\deg(x-\phi)$$
where we view $x$ and $\phi$ as elements of 
$\Hom(A,\hat{A})^+\ot\C$. 
Let $\phi_0:A\to\hat{A}$ be a polarization. Then
the map $x\mapsto \phi_0^{-1}\circ x$ gives an isomorphism
$\Hom(A,\hat{A})^+\ot\Q\simeq \End(A)^+\ot\Q$ (and the corresponding isomorphism
of $\C$-vector spaces). Furthermore, this isomorphism rescales $\deg$ by the
constant $\deg(\phi_0)$. It remains to apply Proposition \ref{deg-inj-prop}.
\ed

\begin{rem} When the ground field is $\C$ we
can identify $\NN(A)\ot\Q$ with the subspace of algebraic cycles in $H^*(A,\Q)$ via the Chern
character and $\NS(A)\ot\Q$ with algebraic cycles in $H^2(A,\Q)$. Then $\ell$ is induced by the
exponential map $H^2(A,\Q)\to H^*(A,\Q)$.
\end{rem}

\subsection{Characterization of the index of a line bundle}

Recall that if $L$ is a nondegenerate line bundle on $A$ (i.e., the corresponding map
$\phi_L:A\to\hat{A}$ is an isogeny) then its {\it index} $i(L)$ is defined by the condition
$H^i(A,L)=0$ for $i\neq i(L)$. We will use the following recipe for computing $i(L)$:
it is the number of positive roots of the polynomial $P(n)=\chi(L\ot L_0^n)$,
where $L_0$ is an ample line bundle on $A$
(see \cite[ch.\ 16]{Mum-ab}). The index function $i(\cdot)$ extends uniquely to a
$\Q_{>0}$-invariant function on $\NS(A)_\Q$.

Let $D_A\sub\NS(A)\ot\C$ be the complexified ample cone.
Note that the function $\deg$ and hence $\chi$ does not vanish on $D_A$ (see \cite[Lem.\ A.3]{GLO}).
Since $D_A$ is simply connected, there is a unique continuous branch of the argument
$\Arg(\chi(x))$ on $D_A$, such that for $x=iH$, where $H$ is an ample class (an element of the ample
cone) we have
$\Arg(\chi(iH))=g\pi/2$, where $g=\dim A$. 
It is easy to see that this branch does not depend on a choice of $H$.
Then for class $x\in\NS(A)\ot\R$ with $\chi(x)\neq 0$ we can define by continuity
the argument $\Arg(\chi(x))$, i.e., we set
$$\Arg(\chi(x))=\lim_{t\to 0+} \Arg(\chi(x+itH)),$$
where $H$ is an ample class. Note that since $\chi(x)$ is real, the number $\Arg(\chi(x))/\pi$ 
is an integer.

\begin{thm}\label{index-thm} For the continuous branch of $\Arg(\chi(\cdot))$ on $D_A$,
satisfying $\Arg(\chi(iH))=g\pi/2$ (where $H$ is ample), one has
$$\Arg(\chi(x))=i(x)\pi$$
for every $x\in\NS(A)\ot\Q$ with $\chi(x)\neq 0$.
\end{thm}

\Pf . First, let us consider the case when $x$ is in the ample cone. For $z\in\C$ we
have $\chi(zx)=z^g\cdot\chi(x)$. Thus, varying $z$ on a unit circle from $1$ to $i$ we obtain
$$\Arg(\chi(ix))=\Arg(\chi(x))+\frac{g\pi}{2}.$$
Since $\Arg(\chi(ix))=g\pi/2$, we obtain that $\Arg(\chi(x))=0$.
Next, assume $x\in\NS(A)\sub\NS(A)\ot\Q$. Then for any ample class $H$ the polynomial
$$P(t)=\chi(x+tH)$$
has $i(x)$ positive roots, counted with multiplicity (see \cite[ch.\ 16]{Mum-ab}).
Let $0<t_1<\ldots<t_r$ be all the positive roots of $P(t)$.
For $t\gg 0$ the class $x+tH$ is ample and so $\Arg\chi(x+tH)=0$.
Now we are going to decrease $t$ until it reaches zero and observe the change of
$\Arg(P(t))=\Arg(\chi(x+tH))$. Note that it can only change when $t$ passes one of the roots $t_j$.
If $t_j$ is a root of $P(t)$ of multiplicity $m_j$, then
for sufficiently small $\eps>0$ one has
$$\Arg(P(t_j-\eps))=\Arg(P(t_j+\eps))+m_j\pi.$$
Adding up the changes we get
$$\Arg(\chi(x))=\Arg(P(0))=i(x)\pi.$$
Since $i(x)$ does not change upon rescaling by a positive rational number,
the assertion for any $x\in\NS(A)\ot\Q$ follows.
\ed

\begin{cor}\label{index-cor} 
For the branch of $\Arg(\deg(\cdot))$ on $D_A$ normalized by $\Arg(\deg(iH))=g\pi$
one has
$$\Arg(\deg(x))=i(x)\cdot 2\pi$$
for any $x\in\NS(A)\ot\Q$ such that $\deg(x)\neq 0$. 
\end{cor}


We will also need some information on the restriction of $\Arg(\chi(\cdot))$ to lines of the form
$iH+\R x\sub D_A$.

\begin{lem}\label{Arg-ineq-lem} 
(i) For any ample class $H\in\NS(A)\ot\Q$ and any $x\in\NS^0(A,\Q)$ let us choose any continuous
branch of $t\mapsto \Arg(\chi(iH+tx))$, where $t\in\R$. 
Then for $0\le t_1<t_2$ one has 
\begin{equation}\label{Arg-ineq}
\Arg(\chi(iH+t_1x))-(g-i(x))\frac{\pi}{2}<\Arg(\chi(iH+t_2x))<\Arg(\chi(iH+t_1x))+i(x)\frac{\pi}{2}.
\end{equation}

\noindent
(ii) For any continuous branch of $\Arg(\deg(\cdot))$ on $D_A$ one has
$$\Arg(\deg(\om))\le \Arg(\deg(iH))+g\pi$$
for any $\om\in D_A$, where $H$ is an ample class.
\end{lem}

\Pf . (i) Indeed, the polynomial
$$P(t)=\chi(iH+tx)=i^g\chi(H+\frac{t}{i}x)$$ 
has all roots purely imaginary, and exactly $i(x)$ of them
in the upper half-plane, counted with multiplicity (see \cite[ch.\ 16]{Mum-ab}). 
Let us write $P(t)=c\cdot (t-z_1)\cdot\ldots\cdot (t-z_g)$.
Since $P(t)\neq 0$ for all $t\in\R$, we can choose for every $j=1\ldots,g$ a continuous branch
of $t\mapsto \Arg(t-z_j)$ along the real line and use the branch
$$\Arg P(t)=\Arg(c)+\Arg(t-z_1)+\ldots+\Arg(t-z_g).$$
Suppose the roots $z_1,\ldots,z_{i(x)}$ are in the upper half-plane while
$z_j$ for $j>i(x)$ are in the lower half-plane. Then
for each $j> i(x)$ the function $t\mapsto \Arg(t-z_j)$
is strictly decreasing and we have
$$\Arg(t_1-z_j)-\frac{\pi}{2}<\Arg(t_2-z_j)<\Arg(t_1-z_j).$$
On the other hand, for $j\le i(x)$ we have
$$\Arg(t_1-z_j)<\Arg(t_2-z_j)<\Arg(t_1-z_j)+\frac{\pi}{2}.$$
Summing up over all the roots gives \eqref{Arg-ineq}.

\noindent
(ii) Applying \eqref{Arg-ineq} to $t_1=0$ and $t_2=1$ we get
$$\Arg(\chi(iH+x))\le\Arg(\chi(iH))+i(x)\frac{\pi}{2}\le \Arg(\chi(iH))+g\frac{\pi}{2}.$$
Since $\deg=\chi^2$ on $\NS$, we get the required inequality for points in $D_A$ with
rational real and imaginary part. The general case follows by continuity.
\ed

\subsection{The group $\BU_{A\times\hat{A}}$}\label{group-sec}

Recall (see \cite{Mukai-spin}, \cite{P-thesis}, \cite{Orlov-ab}, \cite{GLO})
that with every abelian variety $A$ one can associate an algebraic group 
$\BU=\BU_{X_A}$ over $\Q$, where $X_A:=A\times\hat{A}$, as follows.
For every $F\subset\Q$ we define the group of $F$-points $\BU(F)$
as a subgroup of invertible elements of the algebra $\End(X_A)\ot F$ consisting
of 
$$g=\left(\begin{matrix} a & b \\ c & d\end{matrix}\right)\in\End(X_A)\ot F \ \ \text{with }
a\in\Hom(A,A)\ot F, b\in\Hom(\hat{A},A)\ot F, \ \text{etc.},$$
such that 
$$g^{-1}= 
\left(\begin{matrix} \hat{d} & -\hat{b} \\ -\hat{c} & \hat{a}\end{matrix}\right)\in\End(A\times\hat{A})\ot F.$$
The arithmetic subgroup
$$\BU(\Z):=\BU(\Q)\cap\End(A\times\hat{A})$$
is closely related to the group of autoequivalences of $D^b(A)$ (see \cite{Orlov-ab}).
When we view the matrix element $b$ above as a function on $\BU(F)$ we denote it by $b(g)$.

Our point of view is to consider $X_A$ as a ``symplectic object" in the category of abelian varieties
using the skew-symmetric self-duality $\eta_A:X_A\wt{\to}\hat{X_A}$ associated with the biextension
$p_{14}^*\PP\ot p_{23}^*\PP^{-1}$ of $X_A\times X_A$ (see \cite{P-sym}, \cite{P-LIF}). 
Then elements of $\BU(\Z)$
are precisely {\it symplectic automorphisms} of $X_A$, i.e., automorphisms compatible with $\eta_A$.
The development of this point of view in \cite{P-LIF} was to view elements of $\BU(\Q)$
as {\it Lagrangian correspondences} from $X_A$ to itself, which allowed us to define
endofunctors of $D^b(A)$ associated with elements of $\BU(\Q)$ (see \cite[Sec.\ 3]{P-LIF} and
Sec.\ \ref{LI-obj-sec} below).

Note that we have the algebraic 
subgroup $\BT\simeq (\End(A)_\Q)^*\sub \BU$ consisting of diagonal
matrices of the form
$$\left(\begin{matrix} a^{-1} & 0\\ 0 & \hat{a}\end{matrix}\right).$$
 The following facts about the group $\BU$ follow easily from Albert's classification of the
endomorphism algebras of simple abelian varieties (see \cite{P-thesis}, \cite{GLO}).

\begin{lem}\label{alg-group-lem}
(i) Let us fix a polarization $\phi:A\to\hat{A}$ and let
$\BZ\sub 
\BT$ be the algebraic subgroup corresponding to $a\in (\End(A)_\Q)^*$ such that $a$ lies in the
center of $\End(A)_\Q$ and $a^{-1}=\phi^{-1}\hat{a}\phi$. Then the group $\BU$ is an almost direct
product of the semisimple commutant subgroup $S\BU$ and of $\BZ$.

\noindent
(ii) The algebraic group $\BU$ is connected, and the Lie group $\BU(\R)$ is connected (with
respect to the classical topology).
\end{lem}

We denote by $\BU^0\sub \BU$ the Zariski open subset given by the inequality
$\deg(b(g))\neq 0$. 
Note that for any $g\in \BU^0(\R)$ we have $\deg(b(g))>0$ (since the function $\deg$ is
nonnegative on $\Hom(A,\hat{A})\ot\R$).

The following condition on a subset of a group was introduced
in \cite[IV.42]{Weil}  (the term is due to D.~Kazhdan).

\begin{defi}\label{big-defi}
Let $G$ be a group.
A subset $B\sub G$ is called {\it  big}  if  for  any  $g_1, g_2,
g_3\in G$ one has 
$$B^{-1}\cap Bg_1\cap Bg_2\cap Bg_3\neq\emptyset.$$
\end{defi}

This notion is useful because of the following result
(part (i) is due to Weil and part (ii) is a more precise version of \cite[Lem.\ 4.2]{P-maslov}).
   
\begin{lem}\label{big-lem} 
(i) Let $B\sub G$ be a big subset. Then $G$ is isomorphic  to
the abstract group generated by elements $[b]$ for $b\in B$
subject to the relations $[b_1][b_2]=[b_1b_2]$ whenever $b_1b_2\in B$.

\noindent
(ii) Let $Z$ be an abelian group (with the trivial $G$-action). Let
$c,c':G\times G\ra Z$ be a pair of $2$-cocycles such that
$$c(b_1,b_2)=c'(b_1,b_2)$$
for any $b_1,b_2\in B$ with $b_1b_2\in B$.
Let $p:G_c\to G$ (resp. $p':G_{c'}\to G$) be the extension of $G$ by $Z$ 
associated with $c$ (resp., $c'$), and let $\si:G\to G_c$ (resp., $\si':G\to G_{c'}$)
be the natural set-theoretic sections.
Then there is a unique isomorphism
of extensions $i:G_c\to G_{c'}$ such that $i(\si(b))=\si'(b)$ (and identical on $Z$).
\end{lem}
   
\Pf . (i) This is \cite[IV.42, Lem.\ 6]{Weil}.

\noindent
(ii) Note that the subset $p^{-1}(B)\sub G_c$ (resp., $(p')^{-1}(B)\sub G_{c'}$) is big.
Thus, we can define a homomorphism $G_c\to G_{c'}$ be requiring
that it sends $z\si(b)$ to $z\si'(b)$ for $b\in B$, provided we check the compatibility
with the relations 
$$\si(b_1)\si(b_2)=c(b_1,b_2)\si(b_1b_2),$$
$$\si'(b_1)\si'(b_2)=c'(b_1,b_2)\si'(b_1b_2),$$
whenever $b_1, b_2, b_1b_2\in B$. But this boils down to the equality
$c(b_1,b_2)=c'(b_1,b_2)$.
\ed




Next, we will show that the subset $\BU^0(\Q)\sub \BU(\Q)$
(resp., $\BU^0(\R)\sub \BU(\R)$) is big. Note that the subset $\BU^0(\Q)\cap \BU(\Z)$ in
the arithmetic group $\BU(\Z)$ is also big (see Remark \ref{big-Z-rem}).

\begin{lem}\label{dense-lem} 
For any field extension $\Q\sub F$ the set $\BU(F)$ is Zariski-dense in $\BU$.
Hence, the subset $\BU^0(F)\sub \BU(F)$ is big in $\BU(F)$.
\end{lem}

\Pf . Since $\BU$ is connected, density of
$\BU(F)$ follows from \cite[Cor.\ 18.3]{Borel-LAG}.
Thus, for any $g_1,g_2,g_3\in \BU(F)$ the intersection 
$\BU^0\cap \BU^0g_1\cap \BU^0g_2\cap \BU^0g_3$ contains a point of $\BU(F)$.
\ed

The group $\BU$ has two natural parabolic subgroups:
$\BP^+$ is the intersection of $\BU$ with the subgroup of upper-triangular $2\times 2$-matrices
in $\End(A\times\hat{A})_\Q$, and $\BP^-$ is the intersection with the subgroup of lower-triangular
matrices. We also denote by $\BN^+\sub\BP^+$ (resp. $\BN^-\sub\BP^-$) the subgroup
of strictly upper-triangular (resp. strictly lower-triangular) matrices. Note that both $\BN^+$ and
$\BN^-$ are isomorphic to $\NS(A)_{\Q}$.

\begin{lem}\label{P-hom-lem}
Any normal subgroup of $\BU(\Q)$ containing $\BP^-(\Q)$ is the entire $\BU(\Q)$.
\end{lem}

\Pf . Since $\BP^+(\Q)$ is conjugate to $\BP^-(\Q)$ by an element 
\begin{equation}\label{w-phi-eq}
w_\phi=\left(\begin{matrix} 0 & \phi^{-1} \\ -\phi & 0\end{matrix}\right),
\end{equation}
where $\phi:A\to\hat{A}$ is a polarization, it is enough to check that
$\BU(\Q)$ is generated by the subgroups $\BP^+(\Q)$ and $\BP^-(\Q)$.
We can write any $g=\left(\begin{matrix} a & b \\ c & d\end{matrix}\right)\in \BU(\Q)$ with
invertible $a$ as
$$g=\left(\begin{matrix} 1 & 0 \\ ca^{-1} & 1\end{matrix}\right)\cdot
\left(\begin{matrix} a & 0 \\ 0 & \hat{a}^{-1} \end{matrix}\right)\cdot
\left(\begin{matrix} 1 & a^{-1}b \\ 0 & 1\end{matrix}\right).$$
Finally, any element of $\BU^0(\Q)$ has form $gw_\phi$ with $g$ as above.
Thus, the statement follows from Lemma \ref{dense-lem}.
\ed

\subsection{Action of $\BU(\Q)$ on Lagrangian subvarieties}\label{Lag-sec}

Recall that an abelian subvariety $L\sub X_A=A\times\hat{A}$ is {\it isotropic}
if the composition 
$$L\to X_A\rTo{\eta_A}\hat{X_A}\to\hat{L}$$
is zero, where $\eta_A$ is the standard skew-symmetric self-duality. If in addition
$\dim L=\dim A$ then $L$ is called {\it Lagrangian}
(for other equivalent definitions see \cite[Sec.\ 2.2]{P-LIF}).
In this case $\eta_A$ induces an isomorphism $X_A/L\simeq\hat{L}$.

To enumerate all Lagrangian abelian subvarieties in $X_A$ it is convenient to work
in the semisimple category $\Ab_\Q$ of abelian varieties up to isogeny. Note that
abelian subvarieties of $X_A$ are in natural bijection with subobjects of $X_A$ in the category
$\Ab_\Q$. Thus, we can use a similar notion of a Lagrangian subvariety in $\Ab_\Q$.
Now if $L\sub X_A$ is Lagrangian then we have an isomorphism 
$X_A\simeq L\oplus\hat{L}$ in $\Ab_\Q$, which implies that $L$ is isomorphic to $A$ in $\Ab_\Q$.
Thus, we can describe a Lagrangian subvariety (in the category $\Ab_\Q$) as an image of
a morphism $A\to X_A$, i.e., by a pair $(x,y)$, where $x\in\End(A)\ot\Q$, $y\in\Hom(A,\hat{A})\ot\Q$.
The isotropy condition is the equation
$$\hat{y}x=\hat{x}y.$$
The existence of a splitting $X_A\to A$ in $\Ab_\Q$ is equivalent to the condition
$$(\star)\ \ \ (\End(A)\ot\Q) x + (\Hom(\hat{A},A)\ot\Q) y=\End(A)\ot\Q.$$
The pairs $(x_1,y_1)$ and $(x_2,y_2)$ define the same subvariety if and only if
there exists an automorphism $\a$ of $A$ in $\Ab_\Q$ such that $x_2=x_1\a$, $y_2=y_1\a$.
Thus, we obtain an identification of the set of Lagrangian subvarieties in $X_A$ with the set
\begin{equation}\label{LG-set-eq}
\LG(\Q):=\{(x,y) \ |\ \hat{y}x=\hat{x}y, (\star)\}/(x,y)\sim (x\a,y\a),
\end{equation}
where $x\in\End(A)\ot\Q$, $y\in\Hom(A,\hat{A})\ot\Q$ and $\a\in(\End(A)\ot\Q)^*$.
We denote by $(x:y)\in\LG(\Q)$ the equivalence class of 
$(x,y)\in\End(A)\ot\Q\oplus\Hom(A,\hat{A})\ot\Q$.

Fixing a polarization on $A$ we can identify $A$ with $\hat{A}$, so that the dualization gets
replaced by the Rosati involution $x\mapsto x'$ on $\AA:=\End(A)\ot\Q$. 
We claim that the set 
$\LG(\Q)$ can be identified with the set of $\Q$-points of a certain homogeneous projective
variety $\LG$
for the group $\BU$ (a subvariety in the Grassmannian of right rank-1 $\AA$-submodules in
$\AA^2$). Here the action of $\BU$ on $\LG$ is induced by the natural action of
$\End(X_A)_\Q$ on pairs $(x,y)$ (viewed as column vectors).
Consider the point $(0:\phi_0)\in\LG(\Q)$, where $\phi_0:A\to\hat{A}$ is a polarization, 
(the corresponding Lagrangian is $0\times \hat{A}\sub X_A$).
Note that the stabilizer subgroup of is the subgroup $\BP^-\sub\BU$
of lower triangular matrices. Thus, we define
$$\LG=\LG_A=\BU/\BP^-.$$
The fact that the set \eqref{LG-set-eq} is indeed the set of $\Q$-points of $\LG$ follows from the 
transitivity of the action of $\BU(\Q)$ on the set of Lagrangian subvarieties that we will prove below
(see Prop.\ \ref{Lag-transitive-prop}).

We start with the following useful result.

\begin{prop}\label{Lag-auteq-prop} 
For any collection of Lagrangian subvarieties $L_1,\ldots,L_r\sub X_A$
there exists an element $g\in \BU(\Z)$ such that all the Lagrangians
$gL_1,\ldots,gL_r$ are transversal to $\{0\}\times\hat{A}$.
\end{prop}

\Pf . We use an argument similar to the first part of the proof of \cite[Thm.\ 3.2.11]{P-LIF}.
Consider elements in $\BU(\Z)$ of the form $g^+_{nb}$ for some polarization $b:\hat{A}\to A$,
where $n\in\Z$.
Then the condition that $g^+_{nb}L_i$ is transversal to $\{0\}\times\hat{A}$ is equivalent
to $L_i$ being transversal to $g^+_{-nb}(\{0\}\times\hat{A})=\Ga(-nb)$. By \cite[Lem.\ 2.2.7(ii)]{P-LIF},
the latter transversality holds for all $n$ except for a finite number.
\ed

\begin{rem}\label{big-Z-rem}
The above Proposition immediately implies that subset $\BU^0\cap \BU(\Z)$ of the group $\BU(\Z)$ is big (see Sec.\  \ref{group-sec}). 
Indeed, for any given $g_1,\ldots,g_r\in \BU(\Z)$ consider the Lagrangian subvarieties
$L_i=g_i(\{0\}\times\hat{A})\sub X_A$, $i=1,\ldots,n$. Then we can find
$g\in \BU(\Z)$ such that $gL_i=gg_i(\{0\}\times\hat{A})$ for $i=1,\ldots,n$ are transversal to 
$\{0\}\times\hat{A}$. Thus, we get $gg_i\in \BU^0$ as required. The same proof works for
any finite index subgroup $\Ga\sub\BU(\Z)$.
The fact that $\BU^0\cap \Ga$ is a big subset of $\Ga$ was stated in \cite[Lem.\ 4.3]{P-maslov}.
However, the proof in {\it loc. cit.} was not correct: it relied on the absence of compact factors in
$S\BU(\R)$, which is not always the case (see \cite[Cor.\ 5.3.3]{GLO}).
\end{rem}

Lagrangian subvarieties in $X_A$, transversal to $0\times\hat{A}$,
are all graphs $\Ga(f)$ of symmetric homomorphisms $f\in\Hom(A,\hat{A})^+\ot\Q$
(see \cite[Ex.\ 2.2.4]{P-LIF}). This corresponds to points of $\LG(\Q)$ of the
form $(1:f)$, which are precisely $\Q$-points of a Zariski open subset
\begin{equation}\label{NS-LG-subset}
\NS(A)_\Q\simeq \BN^-w_\phi\BP^-/\BP^-\sub \LG,
\end{equation}
where $w_\phi$ is given by \eqref{w-phi-eq}, and $\BN^-\sub\BU$ 
is the subgroup of strictly lower triangular matrices. In other words, the subset \eqref{NS-LG-subset}
is just the $\BN^-$-orbit of the point $(1:0)\in\LG$.

\begin{prop}\label{Lag-transitive-prop} 
(i) The action of $\BU(\Q)$ on the set of Lagrangian subvarieties in $X_A$ is transitive.

\noindent
(ii) The action of $\BU(\R)$ on $\LG(\R)$ is transitive.
\end{prop}

\Pf . (i) The subgroup $\BN^-(\Q)\simeq\NS(A)\ot\Q$ acts on the subset $\NS(A)_\Q\sub\LG$ 
by translations, so the corresponding action on the set of $\Q$-points is transitive. By
Proposition \ref{Lag-auteq-prop}, any point of $\LG(\Q)$ is obtained from a $\Q$-point of
this subset by an action of $\BU(\Z)$, so the required transitivity follows.

(ii) As is well known, it suffices to check triviality of the kernel of the map of Galois cohomology
$H^1(\R,\BP^-)\to H^1(\R,\BU)$. Since $\BP^-$ is a 
semi-direct product of $\prod_i\GL_{n_i}(D_i)$ (where $D_i$ are skew-fields) and of $\G_a^n$, 
in fact, the set $H^1(\R,\BP^-)$ is trivial.
\ed



The description \eqref{LG-set-eq} of $\Q$-points of $\LG$ can be extended to a similar description
of $\LG(F)$, where $F=\R$ or $\C$, so we can still use homogeneous coordinates $(x:y)$,
where $x\in\End(A)\ot F$, $y\in\Hom(A,\hat{A})\ot F$ to describe points of $\LG(F)$.

The complexified ample cone $D_A\sub\NS(A)\ot\C$ is a hermitian
symmetric space (a tube domain) with the group of isometries $\BU(\R)$ 
(see \cite[Sec.\ 5]{Mukai-spin},\cite[Sec.\ 8]{GLO}).
Namely, the group $\BU(\R)$ acts on $D_A$ by
\begin{equation}\label{fract-lin-eq}
g(\om)=(c+d\om)(a+b\om)^{-1}.
\end{equation}
This action is well defined since $\deg(a+b\om)\neq 0$ for $\om\in D_A$ (see
\cite[Lem.\ A3]{GLO}). Furthermore, it is transitive and the stabilizer of a point $\om\in D_A$ is
a maximal compact subgroup of $\BU(\R)$ (see \cite[Thm.\ A1]{GLO}). 
Also, the natural embedding
$$D_A\hra \LG(\C): \om\mapsto (1:\om)$$
is $\BU(\R)$-equivariant.

\section{LI-functors and central extensions}\label{LI-sec}

\subsection{LI-objects and functors}\label{LI-obj-sec}

Recall that every object $K\sub D^b(A\times A)$ gives rise to a functor of Fourier-Mukai type
$$\Phi_K: D^b(A)\to D^b(A): F\mapsto Rp_{2*}(p_1^*F\ot^{{\mathbb L}} K),$$
where $p_1$ and $p_2$ are projections of $A\times A$ to its factors (we refer to $K$ as the
{\it kernel} of the functor $\Phi_K$). The composition
$\Phi_{K_1}\circ\Phi_{K_2}$ corresponds to the convolution of kernels $K_2\circ_A K_1$
(see \cite{Mukai-Fourier}, our notation is as in \cite{P-ker-alg}).

Recall that in \cite{P-LIF} we have extended the relation between autoequivalences of $D^b(A)$ and the group $\BU(\Z)$ (see \cite{P-sym}, \cite{Orlov-ab}) to a construction of endofunctors of $D^b(A)$ (given
by kernels on $A\times A$) associated with elements of $\BU(\Q)$, suitably enhanced. Namely,
with every element $g\in \BU(\Q)$ we associate its graph $L(g)\sub X_A\times X_A$,
which we view as a Lagrangian subvariety in $X_A\times X_A$ with respect to the 
symplectic self-duality $(-\eta_A)\times\eta_A$ (see \cite[Sec.\ 3.1]{P-LIF}). The corresponding kernel
on $A\times A$ is constructed as a generator of the subcategory of $L(g)$-invariants with
respect to the action of $X_A\times X_A$ on $D^b(A\times A)$. 

More precisely, 
every Lagrangian subvariety $L\sub X_A$ can be equipped with a line bundle $\a$
such that we have an isomorphism of line bundles on $L\times L$
\begin{equation}\label{alpha-line-eq}
\a_{l_1+l_2}\ot\a_{l_1}^{-1}\ot\a_{l_2}^{-1}\simeq \PP_{p_A(l_1),p_{\hat{A}}(l_2)},
\end{equation}
where $p_A:L\to A$ and $p_{\hat{A}}:L\to\hat{A}$ are the projections, and $\PP$
is the Poincar\'e bundle on $A\times\hat{A}$. We refer to $(L,\a)$ as {\it Lagrangian pair},
For every such pair $(L,\a)$
there exists a unique up to an isomorphism endosimple coherent sheaf $S_{L,\a}$ on $A$ together
with an isomorphism
\begin{equation}\label{Lag-invariance-eq}
(S_{L,\a})_{x+p_A(l)}\ot\PP_{x,p_{\hat{A}}(l)}\ot\a_l\simeq (S_{L,\a})_x
\end{equation}
on $L\times A$ (where $l\in L$, $x\in A$), satisfying certain natural compatibility condition.
We view this condition as invariance with respect to the lifting of $L$ to
the {\it Heisenberg groupoid} $\bH=\bH_A$, acting on $D^b(A)$ (and on $D^b(A\times S)$ for
any scheme $S$). By definition, $\bH$ is a Picard groupoid extension of $X_A$ by the stack
of line bundles, so its objects over a scheme $S$ are pairs:
a point $(x,\xi)\in X_A(S)$ and a line bundle $\LL$ on $S$. The group operation is determined by
$$(x_1,\xi_1)\cdot (x_2,\xi_2)=\PP_{x_1,\xi_2}\cdot (x_1+x_2,\xi_1+\xi_2).$$
The action of $(x,\xi)\in X_A(S)\sub \bH(S)$ on $D^b(A\times S)$ is given by the functors
\begin{equation}\label{T-x-xi-eq}
F\mapsto T_{(x,\xi)}(F)=\PP_\xi\ot t_x^*F,
\end{equation}
where $\PP_\xi$ is the line bundle on $A\times S$ corresponding to $\xi\in\hat{A}(S)$.
A choice of a line bundle $\a$ satisfying \eqref{alpha-line-eq} 
gives a lifting of $L$ to a subgroup of $\bH$, and the left-hand side of \eqref{Lag-invariance-eq}
is the result of the action of $l\in L$ on $S_{L,\a}$.


\begin{defi}\label{LI-defi}
{\it LI-objects} are cohomologically pure nonzero objects in $D^b(A)$ that can be equipped with
$(L,\a)$-invariance isomorphism \eqref{Lag-invariance-eq} for some $(L,\a)$ as above.
In fact, they are all of the form $S_{L,\a}^{\oplus n}[m]$ for some $(L,\a)$, $n\in\N$ and $m\in\Z$
(see \cite[Thm.\ 2.4.5]{P-LIF}).
Let $\SH^{LI}(A)$ denote the set of isomorphism classes of LI-objects on $A$.
In this work we work mostly with the set $\ov{\SH}^{LI}(A)$ of LI-objects viewed
up to the action of $\bH(k)$, i.e., up to translations and tensoring with line bundles in $\Pic^0(A)$. 
We will refer to this equivalence relation as $\bH$-equivalence.
\end{defi}

We will use the notation $N\cdot F:=F^{\oplus N}$ for an LI-object $F$. This defines
an action of the multiplicative monoid $\N^*$ on $\ov{\SH}^{LI}(A)$.


\begin{prop}\label{LG-map}
There is a well-defined map 
$$\LG(\Q)\to\ov{\SH}^{LI}(A): L\mapsto S(L)$$
sending a Lagrangian subvariety $L\sub X_A$ to the class of the LI-sheaf
$S_{L,\a}$, where $(L,\a)$ is a Lagrangian pair extending $L$. The map
$$\LG(\Q)\times\N^*\times\Z\to \ov{\SH}^{LI}(A):(L,N,n)\mapsto N\cdot S(L)[n]$$
is a bijection of $\N^*\times\Z$-sets.
\end{prop}

\Pf . The fact that $S(L)$ depends only on $L$ follows from \cite[Lem.\ 2.4.2]{P-LIF}.
The second statement follows from
\cite[Thm.\ 2.4.5]{P-LIF} about the structure of the category of $(L,\a)$-invariants in $D^b(A)$
and \cite[Cor.\ 2.4.11]{P-LIF} stating that $L$ can be recoved from $S_{L,\a}$.
\ed

Recall that for an element $\phi\in\NS(A)\ot\Q\simeq\Hom(A,\hat{A})^+\ot\Q$ the graph
$\Ga(\phi)$ is a Lagrangian subvariety of $X_A$. Furthermore, these graphs
are precisely all Lagrangians $L\sub X_A$ such that the projection $L\to A$ is an isogeny.
The sheaf $S_{\Ga(\phi),\a}$ associated with a Lagrangian pair $(\Ga(\phi),\a)$, is a simple semihomogeneous vector bundle with $c_1/\rk=\phi$ (see \cite{Mukai-bun}).
For $\phi\in\NS(A)\ot\Q$ we denote the $\bH$-equivalence class of this bundle by
\begin{equation}\label{V-phi-eq}
V_\phi=S(\Ga(\phi)).
\end{equation}

The above construction of LI-sheaves can be applied
to Lagrangian subvarieties $L\sub X_A\times X_B$ for a pair of abelian varieties $A$ and $B$,
where we use the symplectic self-duality $(-\eta_A)\times\eta_B$ of $X_A\times X_B$.
We refer to the corresponding Lagrangian pairs $(L,\a)$ as {\it Lagrangian correspondences
from $X_A$ to $X_B$}.
The obtained LI-sheaves $S_{L,\a}$ on $A\times B$ can be used as kernels of {\it LI-functors}
$$\Phi_{L,\a}:=\Phi_{S_{L,\a}}:D^b(A)\to D^b(B).$$
The key property of these functors is that we have canonical isomorphisms
\begin{equation}\label{intertwining-eq}
\Phi_{L,\a}\circ T_{p_1(l)}\simeq \a_l\ot T_{p_2(l)}\circ \Phi_{L,\a}
\end{equation}
for $l\in L$, where $p_1,p_2:L\to X_A$ are two projections.
Note that every exact equivalence $D^b(A)\to D^b(B)$ is given by such an LI-functor with
$L$ being the graph of a symplectic isomorphism $X_A\simeq X_B$ (see \cite{Orlov-ab}).

Let $p_{AB}:L\to A\times B$, $p_{A\hat{A}}:L\to A\times\hat{A}$ and $p_{B\hat{B}}:L\to B\times\hat{B}$
be the projections. The line bundle $\a$ can always be chosen in such a way that its restriction
to the connected component of zero in $\ker(p_{AB})$ is trivial. In this case $S_{L,\a}$ is a
direct summand in 
\begin{equation}\label{S-L-for-eq}
p_{AB*}\left(\a^{-1}\ot p^*_{A\hat{A}}\PP^{-1}\ot p^*_{B\hat{B}}\PP\right)
\end{equation}
(see \cite[Lem.\ 3.2.5]{P-LIF}). In the case when $p_{AB}$ is an isogeny the finite group scheme
$\ker(p_{AB})$ has a canonical central extension $H_L$ by $\G_m$ with the underlying
line bundle $\a|_{\ker(p_{AB})}$. Furthermore, $H_L$ is a Heisenberg group scheme and
\eqref{S-L-for-eq} has a natural $H_L$-action, so that
\begin{equation}\label{S-L-bun-eq}
S_{L,\a}=p_{AB*}\left(\a^{-1}\ot p^*_{A\hat{A}}\PP^{-1}\ot p^*_{B\hat{B}}\PP\right)^I,
\end{equation}
for a maximal isotropic subgroup $I\sub\ker(p_{AB})$ lifted to $H_L$. It follows from the theory of
weight one representations of Heisenberg groups that taking $I$-invariants reduces rank by
the factor of $|\ker(p_{AB})|^{1/2}$, so we get
\begin{equation}\label{S-rank-eq}
\rk S_{L,\a}=\deg(p_{AB}:L\to A\times B)^{1/2}.
\end{equation}
In particular, for $B=0$ we get
\begin{equation}\label{rk-V-phi-eq}
\rk V_\phi=\det(p_A:\Ga(\phi)\to A)^{1/2}.
\end{equation}

\begin{ex}\label{tensoring-ex}
The functor of tensoring with a line bundle $L$ on $D^b(A)$ commutes with the action of $\hat{A}$ and
satisfies 
$$L\ot (t_x^*F)\simeq \PP_{-\phi_L(x)}\ot t_x^*(L\ot F).$$
In fact, it is the LI-functor corresponding to
$g_{-\phi_L}=\left(\begin{matrix} \id & 0 \\ -\phi_L & \id\end{matrix}\right)$. More generally, 
for $\phi\in\NS(A)\ot\Q$ the LI-functor corresponding to the element $g_{-\phi}\in\BN^-(\Q)$
is the functor of tensoring with the semihomogeneous vector bundle $V_\phi$ (up to 
$\bH$-equivalence).
\end{ex}

The above construction gives a map
\begin{equation}\label{kernel-map}
\BU(\Q)\to \ov{\SH}^{LI}(A\times A): g\to S(g)=S(L(g)).
\end{equation}
We denote by $\Phi_g\in\Fun(D^b(A),D^b(A))/\bH$ the functor associated with the kernel $S(g)$,
defined up to composing with a functor of the form $T_{(x,\xi)}$, $(x,\xi)\in X_A$ (on either side).
For each $(x,\xi)\in X_A$ we have (noncanonical) isomorphisms
$$\Phi_g\circ T_{N(x,\xi)}\simeq T_{Ng(x,\xi)}\circ\Phi_g,$$
where $N$ is such that $Ng\in \End(X_A)$.

Note that we have a well defined homomorphism induced by $\Phi_g$
$$\rho(g):\NN(A)\to \NN(A).$$
 


\begin{defi}\label{equiv-defi} 
Let $F$ be a cohomologically pure object of $D^b(A)$ and let $G$ be an endosimple
LI-object. We write
$$F\equiv N\cdot G$$
if there exists $n\in\Z$ such that $F[n]$ and $G[n]$ are sheaves and
$F[n]$ has a filtration of length $N$ such that each consecutive quotient is $\bH$-equivalent
to $S(g_1g_2)$. In the case of sheaves on $A\times A$ we will use the same notation for
the relation between the corresponding endofunctors of $D^b(A)$.
\end{defi}

One of the main results of \cite{P-LIF} is the following calculation of the
convolution of kernels (see \cite[Thm.\ 3.3.4]{P-LIF}):
\begin{equation}\label{U-ker-conv-eq}
S(g_2)\circ_A S(g_1)\equiv N(g_1,g_2)\cdot S(g_1g_2)[\la(g_1,g_2)],
\end{equation}
for some $2$-cocycles $N(g_1,g_2)$ and $\la(g_1,g_2)$ of $\BU(\Q)$
with values in $\N^*$ and $\Z$, respectively. 
\footnote{In \cite[Thm.\ 3.3.4]{P-LIF} we made the assumption $\cha(k)=0$ which implies a stronger
statement: the left-hand side of \eqref{U-ker-conv-eq} is a direct sum of objects $\bH$-equivalent
to the right-hand side. It is easy to see that the same argument 
in the positive characteristic case gives a filtration instead of a direct sum.}
Furthermore, we have
\begin{equation}\label{N-for}
N(g_1,g_2)=\frac{q(L(g_1))^{1/2}q(L(g_2))^{1/2}}{q(L(g_1g_2))^{1/2}},
\end{equation}
where 
\begin{equation}\label{q-g-eq}
q(g)=q(L(g))=\deg(p_1:L(g)\to X_A).
\end{equation}
Also, for $g_1,g_2\in \BU^0(\Q)$ such that $g_1g_2\in \BU^0(\Q)$ one has 
\begin{equation}\label{lambda-for}
\la(g_1,g_2)=-i(b(g_1)^{-1}b(g_1g_2)b(g_2)^{-1}).
\end{equation}
Note that in order for the right-hand side to be well-defined
the argument of $i(\cdot)$ should be symmetric.
This indeed follows from the equality 
$$b_1^{-1}(a_1b_2+b_1d_2)b_2^{-1}=b_1^{-1}a_1+d_2b_2^{-1},$$
where we use the usual notation for the entries of $g_1$ and $g_2$.


\begin{defi} 
We denote by $\wt{\BU(\Q)}$ the central extension of $\BU(\Q)$ by $\Z$ 
associated with the $2$-cocycle $\la(\cdot,\cdot)$. Explicitly
$\wt{\BU(\Q)}=\BU(\Q)\times\Z$ with the product 
$$(g_1,n_1)\cdot (g_2,n_2)=(g_1g_2,n_1+n_2+\la(g_1,g_2)).$$
Note that since the subset $\BU^0(\Q)\sub \BU(\Q)$ is big (see Lemma \ref{dense-lem}),
by Lemma \ref{big-lem}(ii), the formula \eqref{lambda-for} determines
the extension $\wt{\BU(\Q)}$ uniquely up to a unique isomorphism.
\end{defi}

Let us denote by $\ov{\SH}^{LI}(A)/\N^*$ the set of equivalence classes with respect to
the equivalence relation generated by $F\sim N\cdot F$ for some $N\in\N^*$.
By \eqref{U-ker-conv-eq}, the map $g\mapsto S(g)\mod\N^*$ defines a homomorphism of monoids
\begin{equation}\label{central-ext-ker-hom}
\wt{\BU(\Q)}\to \ov{\SH}^{LI}(A\times A)^{op}/\N^*,
\end{equation}
and hence a homomorphism of monoids
\begin{equation}\label{central-ext-phi-hom}
\wt{\BU(\Q)}\to \Fun(D^b(A),D^b(A))/(\bH\times \N^*): g\mapsto \Phi_g,
\end{equation}
where on the right we consider functors up to $\bH$-equivalence and up to replacing $\Phi$ with
$N\cdot\Phi=\Phi^{\oplus N}$.

On the level of numerical Grothendieck groups we can eliminate taking quotient by $\N^*$.
Namely, let us set for $g\in \BU(\Q)$ 
\begin{equation}\label{hat-rho-first}
\hat{\rho}(g)=\frac{\rho(g)}{q(g)^{1/2}}:\NN(A)\ot\R\to\NN(A)\ot\R.
\end{equation}
Then from \eqref{U-ker-conv-eq} and \eqref{N-for} we derive that
$$\hat{\rho}(g_1)\hat{\rho}(g_2)=(-1)^{\la(g_1,g_2)}\hat{\rho}(g_1g_2),$$
where $g_1,g_2\in \BU(\Q)$.
Thus, $\hat{\rho}$ defines a homomorphism from 
$\wt{\BU(\Q)}$ to $\GL(\NN(A)\ot\R)$, trivial on the central subgroup 
$2\Z\sub\Z\sub \wt{\BU(\Q)}$.
Note that the quotient $\wt{\BU(\Q)}/2\Z$ is a double cover of $\BU(\Q)$.
Below we will introduce an algebraic structure on this double cover and will show that 
$\hat{\rho}$ is induced by an algebraic homomorphism defined over $\R$
(see Sections \ref{central-ext-sec} and \ref{spin-action-sec}). 


\subsection{Splittings over subgroups}

We are going to define a splitting of the central extension $\wt{\BU(\Q)}\to \BU(\Q)$ over
the parabolic subgroup $\BP^+(\Q)\sub \BU(\Q)$ of lower-triangular matrices (resp.,
over the subgroup $\BP^-(\Q)$ of upper-triangular matrices). Note that
$\BP^+(\Q)$ is a semidirect product of the subgroups of
strictly upper triangular matrices $\BN^+(\Q)\simeq\NS(A)\ot\Q$
and of diagonal matrices $\BT(\Q)\simeq(\End(A)\ot\Q)^*$.

\begin{prop}\label{lifting-prop}
(i) There exist unique liftings of the subgroups $\BN^+(\Q)$ and $\BN^-(\Q)$ to 
$\wt{\BU(\Q)}$. 
The lifting of the element $g^+_\phi=\left(\begin{matrix} 1 & \phi \\ 0 & 1\end{matrix}\right)$, where
$\phi\in\NS^0(A,\Q)$ is given by $(g^+_\phi, i(\phi))\in \wt{\BU(\Q)}$.
The lifting of the element $g^-_\phi=\left(\begin{matrix} 1 & 0 \\ \phi & 1\end{matrix}\right)$, where
$\phi\in\NS^0(A,\Q)$ is given by $(g^-_{\phi},0)\in\wt{\BU(\Q)}$. The corresponding functor 
$\Phi_{g^-_\phi}$ (defined up to $\bH$-equivalence) is given by tensoring with the
semihomogeneous bundle $V_{-\phi}$ (see \eqref{V-phi-eq}).

\noindent
(ii) For $t=t_a=\left(\begin{matrix} a^{-1} & 0\\ 0 & \hat{a}\end{matrix}\right)\in \BT(\Q)$ we have
(up to $\bH$-equivalence)
$$S(t)=\OO_B$$
for some abelian subvariety $B\sub A\times A$ such that the two projections $p,q:B\to A$ are isogenies.
Hence, the functor $\Phi_{t}$ is of the form $q_*p^*$ (up to $\bH$-equivalence).

\noindent
(iii) For any $t\in\BT(\Q)$ and $g\in \BU(\Q)$ one has $\la(t,g)=0$.
\end{prop}

\Pf . (i) Uniqueness of liftings follows from the fact that there are no non-trivial homomorphisms
from a $\Q$-vector space to $\Z$.
Thus, to check the formula for the lifting of $g^+_\phi$ we have to check that
$$S(g^+_{\phi_2})\circ S(g^+_{\phi_1})=S(g^+_{\phi_1+\phi_2})[i(\phi_1+\phi_2)-i(\phi_1)-i(\phi_2),$$
for $\phi_1,\phi_2\in\NS^0(A,\Q)$ such that $\phi_1+\phi_2\in\NS^0(A,\Q)$.
But 
$$\la(g^+_{\phi_1},g^+_{\phi_2})=-i(\phi_1^{-1}(\phi_1+\phi_2)\phi_2^{-1})=-i(\phi_1^{-1}+\phi_2^{-1}),$$
so we are reduced to showing that
$$i(\phi_1^{-1}+\phi_2^{-1})=i(\phi_1)+i(\phi_2)-i(\phi_1+\phi_2).$$
Since 
$$i(\phi_1^{-1}+\phi_2^{-1})=i(\phi_1(\phi_1^{-1}+\phi_2^{-1})\phi_1)=i(\phi_1+\phi_1\phi_2^{-1}\phi_1),$$
this follows from \cite[Prop.\ 15.8]{P-ab-var} (taking into account that $i(-x)=g-i(x)$).

Since the composition of functors $\ot V_{\phi_1}$ and $\ot V_{\phi_2}$ is again tensoring with
a bundle that has a filtration with consecutive quotients 
$\bH$-equivalent to $V_{\phi_1+\phi_2}$, the assertion
about the lifting of $g^-_{\phi}$ follows (cf.\ Ex.\ \ref{tensoring-ex}).

\noindent
(ii) Assume first that $a\in\End(A)$. Then $L(t_a)\simeq A\times\hat{A}$ and its embedding into
$X_A\times X_A$ is given by 
$$(x,\xi)\mapsto (ax, \xi, x, \hat{a}\xi).$$
This implies that $(L(t_a),\OO)$ is a Lagrangian correspondence from $X_A$ to itself, 
so \eqref{S-L-for-eq} in this case gives that
$$S_{L(t_a),\OO}\simeq (a,\id_A)_*\OO_A$$
and the corresponding functor $\Phi_{t_a}$ is the pull-back functor $a^*$.
Similarly, if $a^{-1}\in\End(A)$ then
$$S_{L(t_a),\OO}\simeq (\id_A,a^{-1})_*\OO_A$$
and the corresponding functor $\Phi_{t_a}$ is the push-forward functor $(a^{-1})_*$.
The general case is obtained by combining these two.

\noindent
(iii) We have to check that the convolution $S(g)\circ_A S(t)$ is a sheaf. Indeed, using the form
of $S(t)$ from (ii) we obtain
$$S(g)\circ_A S(t)\simeq (\id_A\times q)_*(\id_A\times p)^*S(g),$$
where $p,q:B\to A$ are isogenies.
\ed

\begin{cor}\label{lifting-cor}
There is a unique splitting of the central extension $\wt{\BU(\Q)}\to \BU(\Q)$ over
$\BP^+(\Q)\sub \BU(\Q)$ (resp., over $\BP^-(\Q)$), which maps $t\in \BT(\Q)$ to $(t,0)\in\wt{\BU(\Q)}$.
\end{cor}

\subsection{Identifying central extensions}\label{central-ext-sec}

Recall that $D_A\sub \NS(A)\ot\C\simeq\Hom(A,\hat{A})^+\ot\C$ denotes the 
complexified ample cone of $A$.

Consider the function $\De:\BU(\R)\to\OO^*(D_A)$ given by 
$$g=\left(\begin{matrix} a & b\\ c& d\end{matrix}\right)\mapsto \De(g)(\om)=\deg(a+b\om),$$
where $\om\in D_A$. 

\begin{lem}\label{De-coc-lem} 
For $g_1,g_2\in \BU(\R)$ one has
\begin{equation}\label{De-cocycle-eq}
\De(g_1g_2)(\om)=\De(g_1)(g_2(\om))\cdot\De(g_2)(\om),
\end{equation}
i.e., $\De$ is a $1$-cocycle.
\end{lem}

\Pf . This follows from the identity
$$a+b\om=(a_1+b_1g_2(\om))(a_2+b_2\om),$$
where $g_i=\left(\begin{matrix} a_i & b_i\\ c_i& d_i\end{matrix}\right)$ for $i=1,2$
and $g_1g_2=\left(\begin{matrix} a & b\\ c& d\end{matrix}\right)$.
\ed

Since $D_A$ is contractible, we have an exact sequence of $\BU(\R)$-modules
$$0\to\Z\to \OO(D_A)\rTo{\exp(2\pi i\cdot ?)}\OO^*(D_A)\to 0.$$
Applying the boundary homomorphism $H^1(\BU(\R))\to H^2(\BU(\Z))$ to the
$1$-cocycle $\De(g)^{-1}$ we obtain a central extension $U^{\De}$ of $\BU(\R)$ by $\Z$.
Explicitly, 
$$U^{\De}=\{(g,f)\in \BU(\R)\times\OO(D_A) \ |\ \De(g)=\exp(-2\pi if)\}.$$
The multiplication rule on $U^{\De}$ uses the cocycle condition on $\De$: we set
$$(g_1,f_1)\cdot (g_2,f_2)=(g_1g_2,f_1(g_2(\cdot))+f_2).$$

\begin{thm}\label{extension-thm} 
There is a homomorphism $\iota:\wt{\BU(\Q)}\to U^{\De}$,
lifting the natural embedding $\BU(\Q)\to \BU(\R)$ and
sending $n\in\Z\sub \wt{\BU(\Q)}$ to $(1,n)\in U^{\De}$. 
This homomorphism is uniquely characterized by the condition that for
$g\in \BU^0(\Q)$ one has
$$\iota(g,0)=(g,f),$$
where
$$\lim_{n\to\infty}\Re f(inH)=-\frac{g}{2}$$
for any ample class $H$.
\end{thm}

\Pf . First, we are going to define a section $\si:\BU^0(\R)\to U^{\De}$ of the projection $U^{\De}\to \BU(\R)$ over
the open subset $\BU^0(\R)\sub \BU(\R)$ consisting of $g$ with $\deg(b(g))\neq 0$.
Note that for $g=\left(\begin{matrix} a & b\\ c & d\end{matrix}\right)\in \BU^0(\R)$ one has
$$\De(g)(\om)=\deg(a+b\om)=\deg(b)\cdot \deg(b^{-1}a+\om).$$
Since $\deg(b)>0$, to define $\si(g)=(g,f^{\si}_g)$ amounts to choosing a branch of the argument for
$\deg(b^{-1}a+\om)^{-1}$. Let us choose the branch of the argument of $\deg(b^{-1}a+\om)$
in such a way that 
$$\lim_{n\to+\infty}\Arg(\deg(b^{-1}a+inH))=\pi\cdot g,$$
where $H$ is an ample class and set $\Arg(\De(g)(\om)^{-1})=-\Arg(\deg(b^{-1}a+\om))$.
Then we set $\iota(g,0)=\si(g)$ for $g\in \BU^0(\Q)$. Since $\BU^0(\Q)$ is big in $\BU(\Q)$, by
Lemma \ref{big-lem}, it remains to show that for $g_1,g_2\in \BU^0(\Q)$ such that
$g_1g_2\in \BU^0(\Q)$ one has
$$\si(g_1)\si(g_2)=\si(g_1g_2)\cdot (1,\la(g_1,g_2)).$$
In other words, we have to check that
$$f^{\si}_{g_1}(g_2(\om))+f^{\si}_{g_2}(\om)=f^{\si}_{g_1g_2}(\om)+\la(g_1,g_2),$$
or equivalently, that with the above choice of $\Arg(\De(g))$ one has
\begin{equation}\label{Arg-identity}
\Arg(\De(g_1)(g_2(\om))+\Arg(\De(g_2)(\om))=\Arg(\De(g_1g_2)(\om))-2\pi\cdot\la(g_1,g_2).
\end{equation}
It is enough to check the equality of the limits of both sides for $\om=inH$ as $n$ goes to infinity
(where $H$ is an ample class).
Let $g_i=\left(\begin{matrix} a_i & b_i\\ c_i & d_i\end{matrix}\right)$ for $i=1,2$.
Note that 
$$\lim_{n\to\infty} g_2(inH)=d_2b_2^{-1}.$$
Thus, \eqref{Arg-identity} reduces to the equality
$$\Arg(\De(g_1)(d_2b_2^{-1}))=-2\pi\la(g_1,g_2)=i(b_1^{-1}b(g_1g_2)b_2^{-1}).$$
But 
\begin{align*}
&\Arg(\De(g_1)(d_2b_2^{-1}))=\Arg(\deg(b_1^{-1}a_1+d_2b_2^{-1})=
\Arg(\deg(b_1^{-1}b(g_1g_2)b_2^{-1})=\\
&2\pi\cdot i(b_1^{-1}b(g_1g_2)b_2^{-1})
\end{align*}
by Corollary \ref{index-cor}.
\ed

The central extension $U^\De\to \BU(\R)$  has a natural continuous splitting over
the subgroup $\BP^-(\R)\sub \BU(\R)$.
Indeed, for $g\in \BP^-(\R)$ we have $\De(g)=\deg(a)>0$, so we can lift $g$ to 
$$\si_{\BP^-}(g)=(g,-\frac{1}{2\pi i}\log(\deg(a))),$$ 
where we choose $\log(\deg(a))$ to be in $\R$. The following result will be useful for us later.

\begin{lem}\label{U-De-splitting-lem}
The restriction of the above lifting homomorphism $\BP^-(\R)\to U^\De$ to $\BP^-(\Q)$
corresponds via $\iota$ to the lifting homomorphism $\BP^-(\Q)\to\wt{\BU(\Q)}$ considered
in Corollary \ref{lifting-cor}. 
\end{lem}

\Pf . By Proposition \ref{lifting-prop}(i), it is enough to check the compatibility of liftings
on $\BT(\Q)$. In view of Proposition \ref{lifting-prop}(iii) this follows from the equality
$$\si_{\BP^-}(t)\si(g)=\si(tg)$$
for any $g\in \BU^0(\Q)$, where $\si:\BU^0(\R)\to U^\De$ is the section used in the proof of
Theorem \ref{extension-thm}.
\ed

Similarly, the extension $U^\De\to \BU(\R)$ has a natural continuous splitting over $\BP^+(\R)$,
which is the same as before over $\BT(\R)$, and over $\BN^+(\R)$ is described as follows.

\begin{lem}\label{N+lifting-lem}
There is a unique splitting of $U^\De\to \BU(\R)$ over $\BN^+(\R)\simeq\NS(\hat{A},\R)$
which is given by the branch of 
$$\Arg\De^{-1}|_{\BN^+(\R)}=\Arg\deg(1+\psi\om)$$ 
that tends to $0$ as $\om\to 0$, where $\psi\in\NS(\hat{A},\R)\simeq\Hom(\hat{A},A)^+_\R$.
\end{lem}

\Pf . It is straightforward to check that this choice of argument gives a lifting. The uniqueness follows
from the fact that there are no nontrivial homomorphisms from a real vector space to $\Z$.
\ed

Let us consider the induced double cover $U^\De/2\Z\to \BU(\R)$.
We are going to introduce an algebraic structure on this group.

\begin{lem}\label{spin-lem} 
Consider a field extension $\Q\sub F$, where either $F=\R$ or $F$ is algebraically closed.
Then for every $g=\left(\begin{matrix} a & b\\ c& d\end{matrix}\right)\in \BU(F)$,
the polynomial $\De(g)(\phi)=\deg(a+b\phi)$ on $\NS(A)(F)$ is a complete square (and is nonzero).
\end{lem} 

\Pf . For $g\in \BU^0$ this follows from the equality 
$$\deg(a+b\phi)=\deg(b)\deg(b^{-1}a+\phi)=\deg(b)\chi(b^{-1}a+\phi)^2$$
and the fact that $\deg(b)\ge 0$ in the case $F=\R$.
Viewing the equation \eqref{De-cocycle-eq} as an identity of rational functions on $\NS(A)$,
we see that if $\De(g_1)$ and $\De(g_2)$ are complete squares then
$\De(g_1g_2)$ is a complete square as a rational function on $\NS(A)$, and hence, as a polynomial.
\ed

\begin{defi} Let $\Pol_{\le g}(\NS(A))$ denote the space of polynomials of degree $\le g$ on $\NS(A)$.
We define a double covering $\Spin=\Spin_{X_A}\to \BU$ of algebraic groups over $\Q$
by setting
$$\Spin=\{(g,f)\in \BU\times\Pol_{\le g}(\NS(A)) \ |\ \De(g)=f^2\}$$
with the group law
$$(g_1,f_1)\cdot (g_2,f_2)=(g_1g_2, f_1(g_2(\cdot))\cdot f_2).$$
Here the rational function $f_1(g_2(\cdot))\cdot f_2$ is actually a polynomial since its square
is $\De(g_1g_2)$.
\end{defi}

Note that by Lemma \ref{spin-lem}, the map $\pi:\Spin(\R)\to \BU(\R)$ is a double covering.
We have a natural isomorphism of groups
\begin{equation}\label{U-De-Spin-isom}
U^{\De}/2\Z\to \Spin(\R): (g,f)\mapsto (g,\exp(-\pi i f)).
\end{equation}
We have two natural subgroups in $\Spin(\R)$:
\begin{equation}\label{BU-spin-eq}
\BU(\Q)^{\spin}=\pi^{-1}(\BU(\Q)),\ \ \BU(\Z)^{\spin}=\pi^{-1}(\BU(\Z)).
\end{equation}

\begin{lem}\label{spin-isom-lem}
Consider the homomorphism
$$\ov{\iota}:\wt{\BU(\Q)}/2\Z\rTo{\sim}\BU(\Q)^{\spin}\sub \Spin(\R)$$
induced by $\iota:\wt{\BU(\Q)}\to U^{\De}$ (see Theorem \ref{extension-thm})
and the isomorphism \eqref{U-De-Spin-isom}.
Then for $g=\left(\begin{matrix} a & b \\ c& d\end{matrix}\right)\in \BU^0(\Q)$ we have
$$\ov{\iota}(g,0)=(g,\sqrt{\deg(b)}\cdot \chi(b^{-1}a+\phi)),$$
with $\sqrt{\deg(b)}>0$.
\end{lem}

\Pf . By Theorem \ref{extension-thm}, $\ov{\iota}(g,0)=(g,f)$,
where $f(\phi)$ is the square root of $\De(g)(\phi)=\deg(b)\cdot \deg(b^{-1}a+\phi)$
with the property
$$\lim_{n\to+\infty}\Arg f(inH)=\frac{\pi\cdot g}{2} \mod 2\pi\Z.$$
Since $\Arg \chi(b^{-1}a+inH)$ has the same limit as $n\to+\infty$, the assertion follows.
\ed

\begin{rems}\label{spin-rem}
1. If for a field extension $\Q\sub F$ there is a multiplicative norm $\Nm$ on $\End(A)\ot F$ 
such that $\Nm^2=\deg$ then the map $g\mapsto (g,\Nm(a+b\om))$ defines a splitting of the
extension $\Spin\to \BU$ over $F$. 
For example, if $A=E^n$, where $E$ is an elliptic curve without complex multiplication, then
$\End(A)=\Mat_n(\Z)$ and $\deg([M]_A)=\det(M)^2$ for a matrix $M\in\Mat_n(\Z)$. Hence, in
this case the norm $\det(\cdot)$ gives a splitting of the spin-covering.

\noindent
2. The group $\BU(\Z)^{\spin}$ is exactly the group 
$USpin(A\times\hat{A})$ defined by Mukai in \cite{Mukai-spin} 
(the same group is denoted by $\Spin(A)$ in \cite{GLO}). 
\end{rems}

Using the isomorphism
$\wt{\BU(\Q)}/2\Z\simeq\BU(\Q)^{\spin}$
we can define a homomorphism
\begin{equation}
\hat{\rho}:\BU(\Q)^{\spin}\to\GL(\NN(A)\ot\R)
\end{equation}
such that $\hat{\rho}(\ov{\iota}(g,0))$ is the operator $\hat{\rho}(g)$ (see \eqref{hat-rho-first}).

\subsection{The action on LI-objects}\label{action-LI-obj-sec}


Recall that with a Lagrangian correspondence from $X_A$ to itself extending
a symplectic isomorphism $g:X_A\to X_A$ in $\Ab_\Q$ we associate an endofunctor
$\Phi_g$ of $D^b(A)$, defined up to $\bH$-equivalence (see Sec.\ \ref{LI-obj-sec}). 
We are going to use these endofunctors to define
an action of $\wt{\BU(\Q)}$ on some extension of $\ov{\SH}^{LI}(A)$ (see Corollary \ref{LI-action-cor}).

\begin{thm}\label{LI-sh-thm} 
(i) For an element $g\in \BU(\Q)$ and a Lagrangian subvariety $L\sub X_A$ we have
\begin{equation}\label{Phi-g-L-eq}
\Phi_g(S(L))\equiv N(g,L)\cdot S(gL)[\la(g,L)]
\end{equation}
for some $\la(g,L)\in\Z$ and $N(g,L)\in\N^*$, where we use Def.\ \ref{equiv-defi}. 

\noindent
(ii) If $L=\Ga(\phi)$ for an isogeny $\phi\in\Hom(A,\hat{A})^+\ot\Q$ and if
$g=\left(\begin{matrix} a & b\\ c& d\end{matrix}\right)$
satisfies $\deg(b)\neq 0$, $\deg(a+b\phi)\neq 0$ and $\deg(c+d\phi)\neq 0$,
then we have
\begin{equation}\label{N-g-L-eq}
N(g,L)=\deg(a+b\phi)^{1/2}\cdot q(g)^{1/2}\cdot\frac{\rk V_{g(\phi)}}{\rk V_\phi},
\end{equation}
where $q(g)$ is given by \eqref{q-g-eq}, $\rk V_\phi$ is given by \eqref{rk-V-phi-eq}, and
\begin{equation}\label{la-g-phi-eq}
\la(g,\Ga(\phi))=-i(b^{-1}a+\phi).
\end{equation}
\end{thm}

\Pf . (i) Let us extend $L$ and $L(g)$ to Lagrangian pairs $(L,\a)$ and $(L(g),\b)$.
By \cite[Thm.\ 3.2.11]{P-LIF}, applied to the Lagrangian correspondence $(L(g),\b)$ and to 
$(L,\a)$ viewed as a Lagrangian correspondence from $0$ to $X_A$, we obtain
$$\Phi_{L(g),\b}(S_{L,\a})=S_{L(g)\circ L,\b\circ\a}[i]$$
for some $i\in\Z$. As in \cite[Thm.\ 3.2.14]{P-LIF} 
one can check that $i$ does not depend on $\a$ and $\b$.
Next, we have to relate the composed Lagrangian correspondence
$S_{L(g)\circ L,\b\circ\a}$ with $S(gL)$. Here we use the definition of the composition
of Lagrangian correspondences from \cite[Sec.\ 3]{P-LIF}. Note that the result is a
{\it generalized Lagrangian correspondence} in the sense of \cite[Def.\ 3.1.1]{P-LIF}
We are going to apply \cite[Prop.\ 2.4.7(ii)]{P-LIF} to the generalized Lagrangian
$Z:=L(g)\circ L\rTo{j} X_A$. 
Note that $Z\sub L(g)\sub X_A\times X_A$ is the preimage of $L$ under the first projection
$p_1:L(g)\to X_A$, and the homomorphism $j:Z\to X_A$ is induced by the second
projection $p_2:L(g)\to X_A$. By   \cite[Prop.\ 2.4.7(ii)]{P-LIF}, we have
$$S_{Z,\b\circ\a}\equiv n^{1/2}\cdot |\pi_0(Z)|^{1/2}\cdot S(j(Z_0))$$
in $\ov{\SH}^{LI}(A\times A)$, 
where $n=|\pi_0(j(Z))|$ (here $Z_0$ is the connected component of $0$ in $Z$).
By definition, we have $j(Z_0)=gL$. Thus, we deduce \eqref{Phi-g-L-eq} with
$$N(g,L)=|\pi_0(Z)|^{1/2}\cdot n^{1/2}.$$
Also, by \cite[(2.4.12)]{P-LIF}, we have $n=\deg(Z_0\to j(Z_0))$.

\noindent(ii)
Now assume that $L=\Ga(\phi)$ and that $g(\phi)$ is defined and is an isogeny.
Note that for sufficiently divisible $N$ we have an isogeny
\begin{equation}\label{Z0-isog-eq}
A\to Z_0: x\mapsto (Nx, N\phi x, N(a+b\phi)x, N(c+d\phi)x)\in L(g)\sub X_A\times X_A.
\end{equation}
In particular, both projections from $Z_0$ to $A$ are isogenies.
Let us consider the commutative diagram of isogenies
\begin{diagram}
Z_0 &\rTo{}& Z\\
\dTo{}&&\dTo{p_{A,2}}\\
j(Z_0)&\rTo{p_A}& A
\end{diagram}
where $p_{A,2}$ is the composition $Z\to L(g)\rTo{p_2} X_A\rTo{p_A} A$.
Considering the degrees we obtain
$$\deg(p_{A,2}:Z\to A)=|\pi_0(Z)|\cdot \deg(p_{A,2}|_{Z_0})=|\pi_0(Z)|\cdot
\deg(j(Z_0)\to A)\cdot n.$$
Recall that $j(Z_0)=gL$, so we get
$$N(g,L)=\frac{\deg(p_{A,2}:Z\to A)^{1/2}}{\deg(gL\to A)^{1/2}}.$$
Now let us consider the projection 
$p_{A,1}:Z\to L(g)\rTo{p_1} X_A\rTo{p_A} A$.
Using the isogeny \eqref{Z0-isog-eq} we see that
$$Np_{A,2}|_{Z_0}=N(a+b\phi)p_{A,1}|_{Z_0}.$$
Hence, 
\begin{equation}\label{ratio-of-projections}
\frac{\deg(p_{A,2}:Z\to A)}{\deg(p_{A,1}:Z\to A)}=\frac{\deg(p_{A,2}|_{Z_0})}{\deg(p_{A,1}|_{Z_0})}=
\deg(a+b\phi).
\end{equation}
Note that $p_{A,1}$ factors through the projection $Z\to L$ and we have a cartesian square
\begin{diagram}
Z &\rTo{}& L(g)\\
\dTo{}&&\dTo{p_1}\\
L&\rTo{}& X_A
\end{diagram}
which shows that $\deg(Z\to L)=\deg(p_1:L(g)\to X_A)=q(L(g))$. Thus,
$$\deg(p_{A,1}:Z\to A)=q(L(g))\cdot\deg(L\to A)$$
and \eqref{ratio-of-projections} can be rewritten as
$$\deg(p_{A,2}:Z\to A)=\deg(a+b\phi)\cdot q(L(g))\cdot\deg(L\to A).$$
Therefore,
$$N(g,L)=\deg(a+b\phi)^{1/2}\cdot q(L(g))^{1/2}\frac{\deg(L\to A)^{1/2}}{\deg(gL\to A)^{1/2}}.$$
Recalling that $\rk S(L)=\deg(L\to A)^{1/2}$ we obtain
\eqref{N-g-L-eq}.

Finally, to compute $\la(g,\Ga(\phi))$ we apply \cite[Prop.\ 3.2.9]{P-LIF}. Namely, we have to consider
the fibered product $\Ga(\phi)\times_A L(g)$ where we use the first projection $L(g)\to A$. Note that we have an isogeny
\begin{equation}\label{la-isog-eq}
A\times\hat{A}\to (\Ga(\phi)\times_A L(g))_0: (x,\xi)\mapsto 
((Nx, N\phi x), (Nx, N\xi, N(ax+b\xi), N(cx+d\xi))),
\end{equation}
where $N$ is sufficiently divisible.
Next, we set 
$$F=\ker(\Ga(\phi)\times_A L(g)\rTo{\ga} A),$$ 
where $\ga$ is induced by the projection to $L(g)$ followed by $L(g)\rTo{p_2} X_A\to A$. 
Note that the composition of $\ga$ with the isogeny \eqref{la-isog-eq} is given
by $(x,\xi)\mapsto N(ax+b\xi).$ Hence, we have an isogeny
\begin{equation}\label{F0-isog-eq}
A\to F_0: x\mapsto ((Nx, N\phi x), (Nx, -Nb^{-1}ax, 0, N(cx-db^{-1}ax))).
\end{equation}
By  \cite[Prop.\ 3.2.9]{P-LIF}, we have
$$\la(g,\Ga(\phi))=-i(g_0\circ f_0^{-1}),$$
where $f_0:F_0\to A$ is the natural projection and for $(l,m)\in F_0\sub \Ga(\phi)\times L(g)$,
$$g_0(l,m)=p_{\hat{A}}(l)-p_{\hat{A},1}(m),$$
where $p_{\hat{A}}:\Ga(\phi)\to \hat{A}$ is the natural projection and $p_{\hat{A},1}$ is the
composition $L(g)\rTo{p_1} X_A\to\hat{A}$. Thus, the compositions of $f_0$ and $g_0$
with the isogeny \eqref{F0-isog-eq} are $x\mapsto Nx$ and $x\mapsto N(\phi+b^{-1}ax)$,
respectively. Hence,
$$g_0\circ f_0^{-1}=\phi+b^{-1}a$$ 
as required.
\ed

Let us set 
$$\ov{\SH}^{LI}(A)_\R=\ov{\SH}^{LI}(A)\times\R_{>0}/\N^*,$$
where $n\in\N^*$ acts by $(F,r)\mapsto (nF,n^{-1}r)$.
Then the bijection of Proposition \ref{LG-map} extends to a bijection of
$\R_{>0}\times\Z$-sets
$$\LG(\Q)\times\R_{>0}\times\Z\rTo{\sim} \ov{\SH}^{LI}(A)_\R.$$

\begin{cor}\label{LI-action-cor} There is an action of $\wt{\BU(\Q)}$ on $\ov{\SH}^{LI}(A)_\R$, commuting
with $\R_{>0}$-action, such that $(g,n)$
acts by 
$$F\mapsto q(L(g))^{-1/2}\cdot\Phi_g(F)[n].$$
For $g_1,g_2\in \BU(\Q)$ and $L\in\LG(\Q)$ we have
$$\la(g_1,g_2(L))+\la(g_2,L)=\la(g_1,g_2)+\la(g_1g_2,L).$$
Also, the maps $\Phi_g$ induce an action of $\wt{\BU(\Q)}$
on $\ov{\SH}^{LI}(A)/\N^*\simeq \LG(\Q)\times\Z$.
\end{cor}

Note that the natural maps 
\begin{equation}\label{SH-LI-NN-Q-map}
\ov{\SH}^{LI}(A)\to \NN(A)\ot\Q/\Q_{>0}
\end{equation} 
\begin{equation}\label{SH-LI-NN-map}
\ov{\SH}^{LI}(A)_\R\to \NN(A)\ot\R
\end{equation} 
associating with an LI-sheaf $F$ its class $[F]$ in $\NN(A)$
are $\wt{\BU(\Q)}$-equivariant, where the action on $\NN(A)\ot\R$ is given by
$\hat{\rho}$ (see \eqref{hat-rho-first}).

\subsection{Action of $\Spin$ on $\NN(A)_\R$}\label{spin-action-sec}

We are going to define an algebraic action of $\Spin$ on $\NN(A)_\R$
inducing the homomorphism $\hat{\rho}$ on $\BU(\Q)^{\spin}\sub\Spin(\R)$.
The idea is to use the algebraicity of the corresponding
projective representation and of the action of an open
subset on a fixed nonzero
vector. We will need the following simple result.

\begin{lem}\label{lifting-lem}
Let $V$ be a vector space over a field $F$, $X$ a scheme (resp., a set),
$\ov{f}:X\to \P(V)$, $\ov{g}:X\to \P(V)$ and $\ov{h}:X\to \P(V)$ be regular morphisms
(resp., maps to the set of $F$-points)
such that the lines $\ov{f}(x)$,$\ov{g}(x)$ and $\ov{h}(x)$ are all distinct and
$\ov{h}(x)\sub\linspan(\ov{f}(x),\ov{g}(x))$ for each $x\in X$. Suppose we have a lifting of $\ov{f}$ to a regular morphism
(resp., map to the set of $F$-points)
$f:X\to V-\{0\}$. Then there exist unique liftings of $\ov{g}$ and $\ov{h}$ to regular morphisms 
(resp., maps to the set of $F$-points)
$g,h: X\to V-\{0\}$ such that $h=f+g$.
\end{lem}

\Pf . Consider a subvariety 
$$Y\sub (V-\{0\})\times (V-\{0\})\times (V-\{0\})$$
consisting of $(v_1,v_2,v_1+v_2)$ such that $v_1$ and $v_2$ are linearly independent,
and a subvariety
$$\ov{Y}\sub (V-\{0\})\times \P(V)\times\P(V)$$
consisting of $(v,L,L')$ such that $v\not\in L$, $v\not\in L'$, $L\neq L'$ and $v\in L+L'$.
Then the natural projection $p:Y\to \ov{Y}$ is an isomorphism.
We have a regular morphism (resp., map to the set of $F$-points) $(f,\ov{g},\ov{h}):X\to\ov{Y}$. Now the components of the corresponding map $X\to Y$ give the required liftings.
\ed

\begin{lem}\label{ell-lem} 
For a symmetric isogeny $\phi\in\NS^0(A,\Q)$ we have
$$\frac{[V_\phi]}{\rk V_\phi}=\ell(\phi)\in\NN(A)\ot\Q,$$
where $V_\phi$ is the semihomogeneous vector bundle \eqref{V-phi-eq} and
$\ell:\NS(A)\ot\Q\to\NN(A)\ot\Q$ is the polynomial map \eqref{ell-map-eq}.
\end{lem}

\Pf . Since $\rk\ell(\phi)=1$,
it suffices to check that the required identity up to proportionality.
Recall that if $(L=\Ga(\phi),\a)$ is a Lagrangian pair 
then the line in $\NN(A)\ot\Q$ corresponding to $V_\phi$ is spanned
by the class of $p_{A*}(\LL)$, where $p_A:L\to A$ is the projection and $\LL=\a^{-1}\ot\PP|_L$
(see \eqref{S-L-bun-eq}). Also, by the definition of a Lagrangian pair, 
$$\La(\a)_{l_1,l_2}\simeq \PP_{p_A(l_1),p_{\hat{A}}(l_2)},$$
so $\phi_{\LL}:L\to\hat{L}$ is given 
$$\phi_{\LL}=\widehat{p_A}\circ p_{\hat{A}}=\widehat{p_{\hat{A}}}\circ p_A,$$
where $p_{\hat{A}}:L\to\hat{A}$ is the projection.
Note that for sufficiently divisible $N$ we have an isogeny
$$i: A\to L: x\mapsto (Nx, N\phi x)$$
and the classes of $p_{A*}(\LL)$ and $[N]_*(i^*\LL)$ in $\NN(A)\ot\Q$ are proportional. 
We have
$$\phi_{i^*(\LL)}=\hat{i}\circ\phi_{\a^{-1}\ot\PP|_L}\circ i=
\widehat{p_A\circ i}\circ p_{\hat{A}\circ i}=N^2\phi.$$
Thus, the class $[i^*(\LL)]\in\NN(A)\ot\Q$ is proportional to $\ell(N^2\phi)=[N]^*\ell(\phi)$.
Hence, the class of $[N]_*(i^*\LL)$ is proportional to
$$[N]_*[N]^*\ell(\phi)=N^{2g}\ell(\phi)$$
as required.
\ed

\begin{thm}\label{alg-spin-action-thm} 
The homomorphism $\hat{\rho}:\BU(\Q)^{\spin}\to\GL(\NN(A)\ot\R)$ 
(see \eqref{hat-rho-first})
extends to an algebraic homomorphism
$$\hat{\rho}:\Spin\to\GL(\NN(A)_\R)$$
defined over $\R$. For $(g,f)\in\Spin(\C)$ and
$\phi\in\NS^0(A,\C)$, such that $g(\phi)$ is defined and belongs
to $\NS^0(A,\C)$, we have
\begin{equation}\label{rho-g-phi-eq}
\hat{\rho}(g,f)(\ell(\phi))=f(\phi)\cdot\ell(g(\phi)).
\end{equation}
\end{thm}

\Pf . First, we observe that Theorem \ref{LI-sh-thm} implies \eqref{rho-g-phi-eq} in the
case when $(g,f)\in \BU(\Q)^{\spin}\sub\Spin(\R)$ with $g\in \BU^0(\Q)$
and $\phi\in\NS^0(A,\Q)$ is such that $g(\phi)$ is defined and belongs to $\NS^0(A,\Q)$.
Indeed, from \eqref{Phi-g-L-eq}, \eqref{la-g-phi-eq}, \eqref{N-g-L-eq} and Lemma
\ref{ell-lem} we obtain in this case
$$\hat{\rho}(g)(\ell(\phi))=(-1)^{i(b^{-1}a+\phi)}|\deg(a+b\phi)|^{1/2}\cdot \ell(g(\phi))=
\deg(b)^{1/2}\cdot \chi(b^{-1}a+\phi)\cdot\ell(g(\phi)).$$
Thus, our claim holds for
$$\ov{\iota}(g,0)=(g,\deg(b)^{1/2}\cdot \chi(b^{-1}a+\phi)).$$
It remains to note that both sides of \eqref{rho-g-phi-eq} change sign when $(g,f)$ gets
multiplied by $-1\in\{\pm 1\}\sub\Spin$.

By \cite[Thm.\ 5.1]{P-maslov}, there is an algebraic homomorphism $\BU\to\PGL(\NN(A)_\Q)$
sending $g\in \BU(\Q)$ to $\rho(g)\mod\Q^*$. Let
us denote the corresponding action of $\BU$ on $\P(\NN(A)_\R)$ by
$$\ov{\kappa}: \BU\times \P(\NN(A)_\R)\to \P(\NN(A)_\R).$$
We also have a map
$$\kappa^\Q: \BU(\Q)^{\spin}\times \NN(A)\ot\R\to\NN(A)\ot\R: (\wt{g},\bv)\mapsto \hat{\rho}(\wt{g})(\bv)
$$
inducing the restriction of $\ov{\kappa}(\R)$ to $\BU(\Q)^{\spin}\times\P(\NN(A)\ot\R)$.
We are going to extend $\kappa^\Q$ to an algebraic morphism using the density
of $\BU(\Q)$ in $\BU$ (see Lemma \ref{dense-lem}).

Note that for a fixed isogeny $\phi$ the right-hand side of \eqref{rho-g-phi-eq}
extends to a regular morphism (defined over $\Q$)
$$\kappa_{\ell(\phi)}: \pi^{-1}(V)\to\NN(A)_\R,$$
where $V\sub \BU^0\sub \BU$ is an open subset of $g\in \BU^0$ 
such that $g(\phi)$ is defined and is an isogeny. Furthermore, as we have seen in the beginning
of the proof, the corresponding map on $\pi^{-1}(V(\Q))$ coincides with the restriction of
$\kappa^\Q$ to $\pi^{-1}(V(\Q))\times\{\ell(\phi)\}$. In particular, the map 
$$\ov{\kappa}_{\ell(\phi)}:\pi^{-1}(V)\to \P(\NN(A)_\R)$$
obtained from $\kappa_{\ell(\phi)}$ is the composition of the projection to $V$
with the restriction of $\ov{\kappa}$ to $\BU\times\{\lan\ell(\phi)\ran\}$
(since we know this on the dense subset $V(\Q)$). 

Now if $\bv\in\NN(A)\ot\Q$ is any vector, linearly independent with $\ell(\phi)$, then 
by Lemma \ref{lifting-lem}, we obtain unique liftings
\begin{equation}\label{two-kappa-eq}
\kappa_\bv, \kappa_{\ell(\phi)+\bv}: \pi^{-1}(V)\to\NN(A)_\R
\end{equation}
of the restrictions of $\ov{\kappa}$ to
$\pi^{-1}(V)\times\{\lan\bv\ran\}$ and $\pi^{-1}(V)\times\{\lan\ell(\phi)+\bv\ran\}$, such that
$$\kappa_{\ell(\phi)+\bv}(\wt{g})=\kappa_{\ell(\phi)}(\wt{g})+\kappa_\bv(\wt{g}).$$
Furthermore, the set-theoretic part of Lemma \ref{lifting-lem}
implies that the maps \eqref{two-kappa-eq} induce the corresponding restrictions of $\kappa^\Q$
on $\pi^{-1}(V(\Q))$. 

Thus, if we consider a basis of $\NN(A)\ot\Q$ of the form $(\ell(\phi),\bv_1,\ldots,\bv_n)$
then combining the maps $\kappa_{\bv_i}$ constructed above we get a regular morphism
$$\hat{\rho}_V: \pi^{-1}(V)\to\GL(\NN(A)_\R)$$
inducing $\hat{\rho}$ on $\pi^{-1}(V(\Q))$. We can cover $\Spin$ with open
subsets of the form $\pi^{-1}(V)\wt{g}$ with $\wt{g}\in\BU(\Q)^{\spin}$ and define a regular
morphism $\pi^{-1}(V)\wt{g}\to\GL(\NN(A)_\R)$ by sending $\wt{h}\wt{g}$ to
$\hat{\rho}_V(\wt{h})\hat{\rho}(\wt{g})$. Using the density of 
$\BU(\Q)^{\spin}$ in $\Spin$, one easily checks that 
these maps glue into the required algebraic homomorphism $\pi^{-1}(V)\to\GL(\NN(A)_\R)$.
\ed

Consider the action of $\Spin(\R)$ on the trivial $\C^*$-bundle $D_A\times\C^*$ over the
domain $D_A$ given by
$$(g,f)\cdot (\om,z)=(g(\om),f(\om)\cdot z),$$
where $(g,f)\in\Spin(\R)$, $\om\in D_A$, $z\in\C^*$.
The map $\ell:D\to\NN(A)\ot\C$ (see \eqref{ell-map-eq}) extends to a $\C^*$-equivariant map
\begin{equation}\label{extended-map-ell}
\ell: D_A\times\C^*\to \NN(A)\ot\C: (\om,z)\mapsto z\cdot\ell(\om).
\end{equation}
From the identity \eqref{rho-g-phi-eq} we immediately get the following result.

\begin{cor}
The map \eqref{extended-map-ell} is $\Spin(\R)$-equivariant.
\end{cor}

\begin{prop}\label{chi-invar-prop} 
For any $x,y\in\NN(A)\ot\C$ and any $\wt{g}\in \Spin(\C)$ one has
\begin{equation}\label{chi-invariance}
\chi(\hat{\rho}(\wt{g})(x),\hat{\rho}(\wt{g})(y))=\chi(x,y).
\end{equation}
\end{prop}

\Pf . Note that the left-hand side of \eqref{chi-invariance} depends only on the image of
$\wt{g}$ in $\BU(\C)$. Let us first consider the case when this image is an element
$g\in \BU(\Q)$. Consider the functor $\Phi=\Phi_{L(g),\a}:D^b(A)\to D^b(A)$ associated with some Lagrangian correspondence $(L(g),\a)$ extending $g$, so that $\Phi$ represents the
$\bH$-equivalence class of $\Phi_g$. Let $\Psi$ be the right adjoint functor to $\Phi$.
By \cite[Prop.\ 3.2.7]{P-LIF}, $\Psi$ differs by a shift from the LI-functor 
associated with some Lagrangian correspondence extending $L(g^{-1})$.
Applying \eqref{U-ker-conv-eq} and \eqref{N-for} for $g_1=g$ and $g_2=g^{-1}$ we obtain
$$\Psi\circ\Phi\equiv N\cdot \Id,$$
where $N=q(g)^{1/2}q(g^{-1})^{1/2}$. 
Since for $F,G\in D^b(A)$ we have an isomorphism
$$\Hom^*(\Phi(F),\Phi(G))=\Hom^*(F,\Psi\Phi(G)),$$
we deduce the equality
\begin{equation}\label{prelim-chi-invar}
\chi(\hat{\rho}(g)([F]),\hat{\rho}(g)([G]))=\frac{q(g^{-1})^{1/2}}{q(g)^{1/2}}\cdot\chi([F],[G]).
\end{equation}
Since $\BU(\Q)$ is dense in $\BU$ (see Lemma \ref{dense-lem}), 
there exists an algebraic character $\varpi:\BU\to\G_m$ 
such that 
$$\chi(\hat{\rho}(g)(x),\hat{\rho}(g)(y))=\varpi(g)\cdot\chi(x,y)$$
for any $g\in \BU(\C)$.
The character $\varpi$ restricts trivially to the semisimple subgroup $S\BU\sub \BU$.
Thus, by Lemma \ref{alg-group-lem}(i), it remains to show the triviality of its restriction to $\BZ$. In fact, we will show directly that
$\varpi(t)=1$ for any $t=\left(\begin{matrix} a^{-1} & 0\\ 0 & \hat{a}\end{matrix}\right)\in \BT(\Q)$,
where $a\in (\End(A)\ot\Q)^*$. Note that this implies that $\varpi|_{\BT}=1$ since $\BT(\Q)$ is dense in 
$\BT$.
It suffices to consider the case when $a\in\End(A)$.
Then the correspondence $L(t)\sub X_A\times X_A$ is the image of the
embedding
$$A\times\hat{A}\to X_A\times X_A: (x,\xi)\mapsto (ax, \xi, x, \hat{a}(\xi)).$$
Hence, in this case $q(t)=\deg(a)$ and $q(t^{-1})=\deg(\hat{a})=\deg(a)$, and our assertion
follows from \eqref{prelim-chi-invar}.
\ed

\begin{cor}\label{chi-ell-cor} 
For $\wt{g}=(g,f_g)\in \Spin(\C)$, $\om\in D_A$ and $x\in\NN(A)\ot\C$ one has
$$\chi(\ell(\om),\hat{\rho}(\wt{g})^{-1}(x))=f_g(\om)\cdot\chi(\ell(g(\om)),x).$$
\end{cor}

\Pf . Indeed, we have
$$\chi(\ell(\om),\hat{\rho}(\wt{g})^{-1}(x))=\chi(\hat{\rho}(\wt{g})(\ell(\om)),x)=f_g(\om)\cdot\chi(\ell(g(\om)),x).$$
\ed

\begin{cor} 
For any $g\in \BU(\Q)$ one has $q(g)=q(g^{-1})$.
\end{cor}


\begin{exs}
1. If $A$ is an abelian variety of dimension $n$ over $\C$ without complex multiplication
then we have $\NS(A)=\Z\cdot H$,
where $H$ is an ample generator, and so $\tau\mapsto \tau\phi_H$ gives an identification
$\HP\to D_A$, where $\HP$ is the upper half-plane. The group $\BU(\R)$ can be identified with 
$\SL(2,\R)$ with the action on $\HP\simeq D_A$ given by fractional-linear transformations 
\eqref{fract-lin-eq}.
Since $\De(g)(\tau\cdot\phi_H)=(a+b\tau)^{2n}$, we have a natural splitting 
$$\SL(2,\R)\to\Spin_{X_A}(\R): g\mapsto (g, (a+b\tau)^n).$$
Furthermore, if $A$ is generic then $\NN(A)\ot\Q$ can be identified with the
$g+1$-dimensional subspace in $H^*(A,\Q)$ spanned by the classes $H^i$, $i=0,\ldots,g$,
and formula \eqref{rho-g-phi-eq} shows that $\SL(2,\R)$ acts on $\NN(A)\ot\R$ as on the standard 
$(g+1)$-dimensional irreducible representation. Assume in addition that $\phi_H$ is a principal
polarization of $A$. Then we can index simple semihomogeneous vector bundles by rational
numbers. Namely, for coprime integers $(r,d)$ with $r>0$ we set
$$V_{r,d}=V_{\frac{d}{r}\phi_H}.$$
From formula \eqref{rk-V-phi-eq} we get in this case $\rk V_{r,d}=r^n$ 
(see also \cite[Rem.\ 7.13]{Mukai-bun}, \cite[ch.\ 12, exer. 2]{P-ab-var}). Hence, by Lemma \ref{ell-lem},
$$\ch(V_{r,d})=\sum_{i=0}^n r^{n-i}d^i\cdot \frac{H^i}{i!}\in H^*(A,\Z).$$
Note that for $r=0$, $d=1$ this formula gives $\ch(\OO_x)$.
Using Hirzebruch-Riemann-Roch formula we get the following
relations for the form $\chi$ on $\NN(A)\ot\Q$:
$$\chi(H^i,H^{n-i})=(-1)^i n!,$$
$$\chi(\ell(\tau\phi_H),[V_{r,d}])=(d-r\tau)^n.$$

\medskip

\noindent
2. Continuing the previous example assume in addition that $n=\dim A=3$ (keeping the assumptions
that $A$ is principally
polarized and generic). 
Then $A$ is the Jacobian of a curve, so $H^2/2$ is an algebraic class.
We claim
that the image of the Chern character $\ch: K_0(A)\to H^*(A,\Q)$ contains the $\Z$-submodule
$K\sub H^*(A,\Q)$ spanned by $(H^i/i!)_{0\le i\le n}$. Indeed, the Chern characters of the
structure sheaves of a point and of the curve span the submodule $\Z H^2/2+\Z H^3/6$. Together
with $\ch(\OO_A)=1$ and $\ch(\OO(H))=\exp(H)$ these classes span the whole $\Z$-submodule $K$. 
On the other hand, using the above
formula we see that for $n\ge 3$
the images of the Chern characters of LI-sheaves (which are all $\bH$-equivalent to either 
$V_{d,r}$ or to $\OO_x$) span a proper $\Z$-submodule in $K$. In particular, the 
LI-objects do not generate $D^b(A)$ in this case.

\medskip

\noindent
3. If $A$ is an elliptic curve over $\C$ with complex multiplication then we have an isomorphism
$$\BU(\R)\simeq\SL(2,\R)\times \BU(1)/\{\pm 1\},$$
where $\{\pm 1\}$ is embedded into the product diagonally. Also, $D_A=\HP$,
the upper half-plane, and $\BU(\R)$ acts on $D_A$ through the projection to $\SL(2,\R)/\{\pm 1\}$.
The spin-covering $\Spin_{X_A}(\R)\to \BU(\R)$ in this case
can be identified with the natural covering
$$\SL(2,\R)\times \BU(1)\to \BU(\R).$$
\end{exs}

\section{Action on stability spaces}\label{stab-sec}

\subsection{Induced $t$-structures and stabilities}\label{t-str-sec}

We refer to \cite{Bridge-stab} for notions related to stability conditions on triangulated categories.
All $t$-structures considered below are assumed to be bounded and nondegenerate
(see \cite{BBD}). All stabilities are assumed to be locally finite and numerical.

We say that a $t$-structure (resp., a slicing or a stability condition) on $D^b(A)$ is $\bH$-invariant, if
it is invariant under any functor $T_{(x,\xi)}$ with $(x,\xi)\in A\times\hat{A}$ 
(see \eqref{T-x-xi-eq}), i.e., under translations and tensoring by $\Pic^0(A)$. 
Note that by \cite[Cor.\ 3.5.2]{P-const}, every {\it full} stability condition is $\bH$-invariant.

The general construction of the induced $t$-structures (resp., stability conditions) 
from \cite{P-const} and \cite{MMS} specializes to the following result on
inducing $\bH$-invariant $t$-structures.

\begin{prop}\label{inducing-prop}
Let $A$ and $B$ be abelian varieties of the same dimension, and let $\Phi:D^b(A)\to D^b(B)$ be
the LI-functor associated with a Lagrangian correspondence $(L,\a)$ from $X_A$ to $X_B$
such that the projections $L\to X_A$ and $L\to X_B$ are surjective, with the right adjoint functor
$\Phi':D^b(B)\to D^b(A)$. Also, let $(D^{\le 0}, D^{\ge 0})$ be an $\bH$-invariant
$t$-structure on $D^b(A)$. Then there is a unique $\bH$-invariant
$t$-structure $(\sideset{^\Phi}{^{\le 0}}{D},\sideset{^\Phi}{^{\ge 0}}{D})$ on $D^b(B)$, such that
\begin{equation}\label{Phi-t-str-incl-eq}
\Phi(D^{[a,b]})\sub \sideset{^\Phi}{^{[a,b]}}{D}
\end{equation}
It is given by
\begin{equation}\label{Phi-t-str-eq}
\sideset{^\Phi}{^{[a,b]}}{D}=\{F\in D^b(B)\ |\ \Phi'(F)\in D^{[a,b]}\}.
\end{equation}

Similarly, if $(P(t))_{t\in\R}$ is a $\bH$-invariant slicing on $D^b(A)$ then there is a unique $\bH$-invariant slicing
$(\sideset{^\Phi}{}{P}(t))_{t\in\R}$ on $D^b(B)$ such that $\Phi(P(t))\sub\sideset{^\Phi}{}{P}(t)$ for 
any $t\in\R$. We have $\sideset{^\Phi}{}{P}(t)=(\Phi')^{-1}(P(t))$.
\end{prop}

\Pf . First, we observe that by Proposition \cite[Prop.\ 3.2.7]{P-LIF},
$\Phi'$ differs by a shift from the LI-functor associated with the transposed correspondence
$(\si(L),\a^{-1})$, where $\si:X_A\times X_B\to X_B\times X_A$ is the permutation of factors.
Hence,
the same argument as in \cite[Lem.\ 3.3.3]{P-LIF} shows that both compositions $\Phi'\circ\Phi$ and $\Phi\circ\Phi'$ are obtained by consecutive extensions from functors $T_{(x,\xi)}$, one of which is the identity functor.

The fact that \eqref{Phi-t-str-eq} defines a $t$-structure
follows from \cite[Thm.\ 2.1.2]{P-const} once we check that in our situation $\Phi'\circ\Phi$ is $t$-exact with respect to the original $t$-structure and $(\Phi\circ\Phi')(F)=0$ implies $F=0$. Indeed, the former
follows from $\bH$-invariance of our $t$-structure. To check the latter property it suffices to consider
the case when $F$ is a coherent sheaf. We observe that the right adjoint functor to $\Phi\circ\Phi'$ sends a structure sheaf of a point $\OO_x$ to a sheaf $K_x$ supported on a finite number of points including $x$. Hence, if $(\Phi\circ\Phi')(F)=0$ then $\Hom(F,K_x)=0$ for all $x\in B$, 
which implies that $F=0$.

The inclusion \eqref{Phi-t-str-incl-eq} follows from the $\bH$-invariance of the original
$t$-structure and from the form of $\Phi'\circ\Phi$. 
The fact that the new $t$-structure is $\bH$-invariant follows from the $\bH$-intertwining property
of LI-functors (see \eqref{intertwining-eq} and \cite[Lem.\ 3.2.4]{P-LIF}).
Now suppose 
$(\sideset{^\Phi}{^{\le 0}}{D_1},\sideset{^\Phi}{^{\ge 0}_1}{D})$ is another $\bH$-invariant
$t$-structure
on $D^b(B)$ such that $\Phi(D^{[a,b]})\sub\sideset{^\Phi}{^{[a,b]}_1}{D}$. Then applying
\cite[Thm.\ 2.1.2]{P-const} again we deduce that
$$D_1^{[a,b]}=\{F\in D^b(A)\ |\ \Phi(F)\in \sideset{^\Phi}{^{[a,b]}_1}{D}\}$$
is a $t$-structure on $D^b(A)$ such that $D^{[a,b]}\sub D_1^{[a,b]}$.
Hence, $D_1^{[a,b]}=D^{[a,b]}$ and we can rewrite \eqref{Phi-t-str-eq} as
$$\sideset{^\Phi}{^{[a,b]}}{D}=\{F\in D^b(B)\ |\ \Phi\Phi'(F)\in \sideset{^\Phi}{^{[a,b]}_1}{D}\},$$
which implies that $\sideset{^\Phi}{^{[a,b]}_1}{D}\sub\sideset{^\Phi}{^{[a,b]}}{D}$,
so these $t$-structures are the same.

The result about slicings is proved analogously.
\ed

Let $\Stab^{\bH}(A)$ denote the space of $\bH$-invariant stability conditions on $A$
(it is known to be nonempty for $\dim A\le 2$).

\begin{defi}\label{stab-action-defi} 
For $g\in \BU(\Q)$ and a stability $\si=(P(\cdot),Z)\in\Stab^{\bH}(A)$ we set
$$g(\si)=(\sideset{^{\Phi_g}}{}{P}(\cdot), Z\circ \hat{\rho}(g)^{-1}),$$
where $\hat{\rho}(g)$ is given by \eqref{hat-rho-first}. 
By Proposition \ref{inducing-prop}, this defines an action of $\wt{\BU(\Q)}$ on $\Stab^{\bH}(A)$,
such that the central element $1\in\Z\sub\wt{\BU(\Q)}$ sends $(P(\cdot),Z)$ to $(P(\cdot)[1],-Z)$.
\end{defi}

The restriction of the above action to the preimage of $\BU(\Z)\sub \BU(\Q)$ is
given by the standard action of the autoequivalence group of $D^b(A)$
on $\Stab(A)$ (see \cite{Bridge-stab}).

\begin{prop} For every $\wt{g}\in\wt{\BU(\Q)}$ the corresponding transformation of
$\Stab^{\bH}(A)$ is an isometry with respect to the generalized metric $d(\cdot,\cdot)$ introduced in
\cite[Prop.\ 8.1]{Bridge-stab}.
\end{prop}

\Pf . Note that the functor $\Phi_g$ sends Harder-Narasimhan constituents
of $E$ with respect to $\si$ to those of $\Phi_g(E)$ with respect to $g(\si)$, and
$Z(\hat{\rho}(g)^{-1}(\Phi_g(E)))$ is a constant multiple (depending only on $g$) of
$Z(E)$. Hence, $d(\si_1,\si_2)\le d(g(\si_1),g(\si_2))$. Applying the same inequality to $g^{-1}$
and the pair $(g(\si_1),g(\si_2))$ we deduce that it is in fact an equality.
\ed

\begin{prop}\label{LI-semistable-prop} 
Any LI-object in $D^b(A)$ is semistable with respect to any full stability. 
\end{prop}

\Pf . Let $E$ be an $(L,\a)$-invariant object in $D^b(A)$, where $(L,\a)$ is a Lagrangian pair
(with $L\sub X_A$), and let $\si=(P(\cdot),Z)$ be a full stability. We can assume that $Z$ takes values in $\Q+i\Q\sub\C$. Indeed, the set of such stabilities is dense in the connected component containing
$\si$, and the semistability of $E$ is a closed condition on $\si$. Then for a dense set of real numbers
$t$ (namely, those with $\tan(\pi t)\in\Q$) the abelian category $P((t,t+1])$ is Noetherian (see
\cite[Prop.\ 5.0.1]{AP}). Applying the construction of \cite{P-const} we obtain
for each such $t$ the associated
{\it constant family} of $t$-structures over any base $S$, which is a certain $t$-structure on
$D^b(A\times S)$,  local over $S$ and such that its heart contains the pull-back of 
$P((t,t+1])$ with respect to the projection $p_1: A\times S\to A$. Let us take as a base $S=L$ and
consider the functor 
$$T_{(L,\a)}: D^b(A)\to D^b(A\times L)$$ 
that associates with $F\in D^b(A)$
the natural family of objects $\FF$ on $L\times A$ such that the
restriction of $\FF$ to $\{l\}\times A$ is $\a_l\ot T_l(F)$ for $l\in L$ ($\FF$ is obtained from $F$
by taking the pull-back with respect to the map $L\times A\to A:(l,x)\mapsto p_A(l)+x$ and
then tensoring the result with a certain line bundle). Since our stability is $\bH$-invariant
(by \cite[Cor.\ 3.5.2]{P-const}), this functor is easily seen to be $t$-exact, i.e., it sends 
$P((t,t+1])$ to the heart of the corresponding constant $t$-structure on $D^b(A\times L)$.
By definition, $(L,\a)$-invariance structure on $E$ is an isomorphism
$$T_{(L,\a)}(E)\simeq p_1^*E.$$
Since both sides are $t$-exact functors of $E$, we deduce that the truncations of $E$ with respect to
our $t$-structure are still $(L,\a)$-invariant. Applying this for an appropriate set of phases $t$
we derive that all Harder-Narasimhan constituents of $E$ are 
$(L,\a)$-invariant. Let $E_0$ be one of them. Suppose $E_0$ has cohomological range $[a,b]$
with respect to the standard $t$-structure. Then $H^bE_0$ and $H^aE_0$ are still
$(L,\a)$-invariant, so we have a nonzero morphism $H^bE_0\to H^aE_0$ (see \cite[Thm.\ 2.4.5]{P-LIF}), which gives rise to a nonzero morphism
$$E_0[b]\to H^bE_0\to H^aE_0\to E_0[a].$$
By semistability of $E_0$ we should have $b\le a$, i.e., $E_0$ is cohomologically pure.
Since $E_0$ is a direct sum of several copies of the generator $S_{L,\a}$, it follows that
$S_{L,\a}$ is also semistable.
\ed

\subsection{$\Z$-covering of $\LG(\R)$}

Recall that the action of $\BU(\R)$ on $\LG(\R)$ is transitive (see Prop.\ \ref{Lag-transitive-prop}), 
so we have an identification
\begin{equation}\label{LG-R-quot-eq}
\LG(\R)\simeq \BU(\R)/\BP^-(\R).
\end{equation}

We have a natural lifting of $\BP^-(\R)$ to a closed subgroup of $U^\De$ (see Lemma
\ref{U-De-splitting-lem}).
Therefore, the homogeneous space $U^\De/\BP^-(\R)$ is a $\Z$-covering of $\LG(\R)$.
Below we will describe this $\Z$-covering explicitly using the homogeneous coordinates
$(x:y)$ on $\LG(\R)$ (see Sec.\ \ref{Lag-sec}).

Namely, with every $L=(x:y)\in\LG(\R)$ we associate a holomorphic function on $D_A$, defined
up to rescaling by a positive constant, 
$$\de(L)(\om)=\deg(\hat{y}-\hat{x}\om)=\deg(\om x-y) \mod \R_{>0},$$
where $\om\in D_A$. Note that if we change $(x:y)$ to $(x\a:y\a)$ then this function gets multiplied by
$\deg(\a)\in\R_{>0}$.
It is easy to see that for $g\in \BU(\R)$ one has
\begin{equation}\label{de-De-eq}
\de(g(0:\phi_0))=\De(g^{-1}) \mod \R_{>0},
\end{equation}
where $\phi_0:A\to\hat{A}$ is a polarization
and $\De$ is the $1$-cocycle of $\BU(\R)$ with values in $\OO^*(D_A)$ defined in
Section \ref{central-ext-sec}.
In particular, $\de(L)(\om)\neq 0$  for all $\om\in D_A$.

\begin{lem}\label{de-LG-lem} 
For $L\in\LG(\R)$ and $g\in \BU(\R)$ one has
$$\de(gL)(g(\om))=\de(L)(\om)\cdot\De(g^{-1})(g(\om))=
\de(L)(\om)\cdot\De(g)(\om)^{-1}.$$
\end{lem}

\Pf . Pick $g'\in \BU(\R)$ such that $L=g'(0:\phi_0)$. Then use \eqref{de-De-eq} and
the cocycle condition for $\De$.
\ed

Note also that if we have a Lagrangian subvariety $L\sub A\times\hat{A}$ then viewing $L$ as a point in
$\LG(\Q)$ we have
\begin{equation}\label{de-L-Lag-eq}
\de(L)(\om)=\deg(\om p_1-p_2)\mod \R_{>0},
\end{equation}
where $p_1:L\to A$ and $p_2:L\to\hat{A}$ are the projections, and we use the polynomial
function $\deg:\Hom(L,\hat{A})\ot\C\to\C$.

\begin{defi} We define the $\Z$-covering
$p:\wt{\LG(\R)}\to \LG(\R)$
by setting
$$\wt{\LG(\R)}=\{(L,f)\in\LG(\R)\times(\OO(D_A)/i\R)\ |\ \de(L)=\exp(2\pi i f) \mod \R_{>0}\}.$$
We also set
$$\wt{\LG(\Q)}:=p^{-1}(\LG(\Q))\sub \wt{\LG(\R)}.$$
When we need to stress the dependence on $A$ we write $\wt{\LG_A(\R)}$ (resp., 
$\wt{\LG_A(\Q)}$).
We have an action of $U^\De$ on $\wt{\LG(\R)}$ given by
$$(g,f_g)\cdot (L,f_L)=(gL, f_L(g^{-1}(\om))+f_g(g^{-1}(\om))).$$
The fact that this action is well defined follows from Lemma \ref{de-LG-lem}.
\end{defi}

\begin{prop}\label{Z-cov-prop} (i) There exists a unique bijection
\begin{equation}\label{Z-cov-bijection}
\ov{\SH}^{LI}(A)/\N^*\to \wt{\LG(\Q)}:F\mapsto \wt{L}_F,
\end{equation}
lifting the natural projection $F\mapsto L_F$ to $\LG(\Q)$,
sending $\OO_x$ to $((0:\phi_0),0)\in\wt{\LG(\Q)}$ (where $\phi_0:A\to\hat{A}$
is a polarization),
and $\wt{\BU(\Q)}$-equivariant, where the action on $\wt{\LG(\Q)}$ is
induced by the embedding $\iota:\wt{\LG(\Q)}\to U^\De$.

\noindent
(ii) Let $V_\phi$ be the semihomogeneous vector bundle associated
with $\phi\in\NS(A)\ot\Q$, so that $L_{V_\phi}=\Ga(\phi)$ (see \eqref{V-phi-eq}). Then
\begin{equation}\label{wt-L-V-phi-eq}
\wt{L}_{V_\phi}=(\Ga(\phi), (2\pi i)^{-1}\cdot\log(\deg(\om-\phi)) \mod i\R),
\end{equation}
where the branch of $\log(\deg(\cdot))$ is normalized by 
$\Im\log(\deg(iH))=\Arg(\deg(iH))=-g\pi$ for ample $H$.
\end{prop}

\Pf . (i) First, let us compute the stabilizer subgroup $\St\sub\wt{\BU(\Q)}$ of the class of 
$\OO_x$ in $\ov{\SH}^{LI}(A)/\N^*$.
By considering the action on the corresponding Lagrangian we see that 
$\St$ is a certain lifting of $\BP^-(\Q)\sub \BU(\Q)$ to $\wt{\BU(\Q)}$. From
the explicit form of the functors $\Phi_t$ for $t\in\BT(\Q)$ (see Prop.\ \ref{lifting-prop}(ii))
we see that these functors preserve $\OO_x$ up to $\bH$-equivalence and $\N^*$.
Therefore, $\St$ is the lifting of $\BP^-(\Q)$ described in Corollary \ref{lifting-cor}.
By Lemma \ref{U-De-splitting-lem}, $\iota(\St)$ is exactly the stabilizer of the point
$((0:\phi_0),0)\in\wt{\LG(\Q)}$. Hence, there is a well-defined
$\wt{LG(\Q)}$-equivariant map \eqref{Z-cov-bijection}. Since this is a map of
$\Z$-torsors over $\LG(\Q)$ (see Prop.\ \ref{LG-map}), it is a bijection.

\noindent
(ii) Assume first that $\phi$ is non-degenerate, i.e., $\phi\in\NS^0(A,\Q)$.
Consider the element $g^+_{\phi^{-1}}\in\BN^+(\Q)$ as in Proposition \ref{lifting-prop}(i). Then
$$g^+_{\phi^{-1}}(0:\phi_0)=(\phi^{-1}\phi_0:\phi_0)=(1:\phi)=\Ga(\phi).$$
By Proposition \ref{lifting-prop}(i), under the canonical lifting of $\BN^+(\Q)$ to $\wt{\BU(\Q)}$ the lifting of $g^+_{\phi^{-1}}$ corresponds to the kernel $S(g^+_{\phi^{-1}})[i(\phi)]$ (note that $i(\phi^{-1})=i(\phi)$). 
On the other hand, its canonical lifting to $U^\De$ is 
$$\wt{g}=(g^+_{\phi^{-1}},-\log(\deg(1+\phi^{-1}\om))/2\pi i),$$
where we use the branch of $\Arg\deg(1+\phi^{-1}\om)$ that tends to $0$ as $\om\to 0$ (see Lemma
\ref{N+lifting-lem}).
By the $\wt{\BU(\Q)}$-equivariance of the map \eqref{Z-cov-bijection}, we obtain that 
the object $V_\phi[i(\phi)]$ is mapped under this map to 
\begin{align*}
&\wt{g}\cdot((0:\phi_0),0)=(\Ga(\phi),\log(\deg(1-\phi^{-1}\om))/2\pi i \mod i\R)=\\
&(\Ga(\phi),(2\pi i)^{-1}\cdot\log(\deg(\om-\phi)) \mod i\R)
\end{align*}
with the same choice of the argument as above.
Recall that if we choose the branch of $\Arg\deg(\om-\phi)$ in such a way
that $\Arg\deg(inH-\phi)$ will be $\pi g$ then we will obtain
$$\Arg\deg(-\phi)=2\pi i(-\phi)=2\pi (g-i(\phi))$$ 
(see Corollary \ref{index-cor}).
Subtracting $2\pi g$ we get the branch that gives the limit
$-2\pi i(\phi)$ as $\om\to 0$ which is exactly what we get for the image of $V_\phi$.
This proves the required statement in the case when $\phi\in\NS^0(A,\Q)$.
The general case follows by using the action of the subgroup $\BN^-(\Q)\sub\wt{\BU(\Q)}$
(see Ex.\ \ref{tensoring-ex}). 
Indeed, this action changes both sides \eqref{wt-L-V-phi-eq} by adding to
$\phi$ an arbritrary element of $\NS(A)\ot\Q$.
\ed

Recall that we can view $\NS(A)\ot\R$ as an open subset of $\LG(\R)$ via the map
$\phi\mapsto (1:\phi)=\Ga(\phi)$. Part (ii) of the above Proposition implies
that we have a commutative diagram
\begin{diagram}
\NS(A)\ot\Q&\rTo{\phi\mapsto V_\phi} & \ov{\SH}^{LI}(A)/\N^*\\
\dTo{}&&\dTo{}\\
\NS(A)\ot\R&\rTo{}& \wt{\LG(\R)}
\end{diagram}
where the right vertical arrow is \eqref{Z-cov-bijection} and the bottom arrow
is the continuous section of the projection $\wt{\LG(\R)}\to\LG(\R)$ over $\NS(A)\ot\R$ given by
\begin{equation}\label{NS-LG-section-eq}
\NS(A)\ot\R\to\wt{\LG(\R)}: \phi\mapsto (\Ga(\phi), f_\phi),
\end{equation}
where $f_\phi\in \OO(D_A) \mod i\R$ is the branch of
$(2\pi i)^{-1}\cdot\log(\deg(\om-\phi)) \mod i\R$ satisfying
$$\lim_{n\to\infty} f_\phi(inH)=-g/2$$ 
for any ample $H$. 

\begin{defi}\label{LG-spin-defi}
We define the double covering $p^{\spin}:\LG^{\spin}(A,\R)\to\LG_A(\R)$ by setting
$\LG^{\spin}(A,\R)=\wt{\LG(\R)}/2\Z$. Explicitly, 
$$\LG^{\spin}(A,\R)=\{(L,\varphi)\in\LG(\R)\times(\OO(D_A)/\R_{>0})\ |\ \de(L)=\varphi^2 \mod \R_{>0}\}.$$
We also set $LG^{\spin}(A,\Q)=(p^{\spin})^{-1}(\LG(\Q))$.
\end{defi}

The isomorphism $\Spin(\R)\simeq U^\De/2\Z$ induces a transitive
action of $\Spin(\R)$ on $\LG^{\spin}(A,\R)$ (and of $\BU(\Q)^{\spin}$ on $\LG^{\spin}(A,\Q)$).

We also have a natural $\wt{\BU(\Q)}$-equivariant
map 
$$\ov{\SH}^{LI}(A)/\N^*\to \NN(A)\ot\Q/\Q_{>0}: F\mapsto [F]\mod\Q_{>0}$$ 
(see \eqref{SH-LI-NN-Q-map}), 
which we can view as a map from $\wt{\LG(\Q)}$ using the bijection \eqref{Z-cov-bijection}.
The equivariance of this map with respect to the $\Z$-action 
implies that it factors through $\LG^{\spin}(A,\Q)$. Furthermore, we claim
that it extends to a continuous $\Spin(\R)$-equivariant map
\begin{equation}\label{LG-spin-NN-map}
\LG^{\spin}(A,\R)\to \NN(A)\ot\R/\R_{>0}
\end{equation}
such that we have a commutative  diagram:
\begin{diagram}
\ov{\SH}^{LI}(A)/\N^*&\rTo{}& \NN(A)\ot\Q/\Q_{>0}\\
\dTo{}&&\dTo{}\\
\LG^{\spin}(A,\R)&\rTo{}& \NN(A)\ot\R/\R_{>0}
\end{diagram}
Indeed, we can define \eqref{LG-spin-NN-map} by sending
$\wt{g}((0:\phi_0),1)$ to $\wt{g}[\OO_x]\mod\R_{>0}$ for $\wt{g}\in\Spin(\R)$.
To check that this map is well defined we observe that $\Q$-points $\BP^-(\Q)$ are dense (with
respect to the classical topology) in the stabilizer $\BP^-(\R)$ of the point $((0:\phi_0),1)\in\LG^{\spin}(A,\R)$.
Since $\BP^-(\Q)\sub\BU(\Q)^{\spin}$ leaves the class $[\OO_x]\in S\NN(A)\ot\R$ 
invariant, this proves our claim.


\begin{lem}\label{NS-section-lem} 
The section \eqref{NS-LG-section-eq} induces a section
\begin{equation}\label{NS-LG-spin-section-eq}
\NS(A)\ot\R\to \LG^{\spin}(A,\R)
\end{equation}
which sends $\phi\in\NS(A)\ot\R$ to $(\Ga(\phi),\chi(\phi-\om)\mod\R_{>0})$.
\end{lem}

\Pf . Since $\chi(\phi-\om)^2=\deg(\phi-\om)=\deg(\om-\phi)$,
this follows from the fact that the argument of $\chi(\phi-inH)$ tends to $-g\pi/2\mod 2\pi\Z$ as 
$n\to\infty$.
\ed

\begin{ex}\label{ell-power-ex}
In the case when $A=E^n$, where $E$ is an elliptic curve without complex multiplication
we can identify $\End(A)$ with the algebra of $n\times n$-matrices over $\Z$, and
$\NS(A)$ with symmetric matrices. Note that for $M\in\End(A)$ we have
$\deg(M)=\det(M)^2$ and for $\phi\in\NS(A)$ we have $\chi(\phi)=\det(\phi)$. In a coordinate-free
notation, if $A=E\ot \La$, where $\La$ is a free $\Z$-module of rank $n$, then elements of $\NS(A)$
can be viewed as $\Z$-valued symmetric bilinear forms on $\La$, and the function $\chi$ is given
by the discriminant.
The group $\BU$ in this case is the symplectic group $\Sp_{2n}$ and the variety
$\LG_A$ is the Lagrangian Grassmannian associated with the $2n$-dimensional symplectic vector
space. Also, $D_A$ is the Siegel 
upper half-plane $\HP_n$ and the covering $U^\De\to\Sp_{2n}$ corresponds to a choice
of argument of $Z\mapsto \det(A+BZ)^2$, where $Z\in\HP_n$ and
$\left(\begin{matrix} A & B\\ C& D\end{matrix}\right)\in\Sp(2n,\R)$.
Thus, $U^\De$ contains the universal covering $\wt{\Sp}(2n,\R)$ of $\Sp(2n,\R)$
as a subgroup of index $2$ (cf. \cite[Ex.\ 4.15]{Orlov-ab}).
Now let us consider our lifting of $\BP^-(\R)$ to $U^\De$. It is easy to check that the restriction of
the projection to $U^\De/\wt{\Sp}(2n,\R)\simeq\{\pm1\}$ to $\GL(n,\R)\sub\BP^-(\R)$ can be identified
with the homomorphism $A\mapsto \sign\det(A)$.
It follows that $U^\De=\wt{\Sp}(2n,\R)\cdot \BP^-(\R)$, and
$\BP^-(\R)\cap\wt{\Sp}(2n,\R)$ is the semidirect product of $\BN^-(\R)$ and of $\GL^+(n,\R)$
(matrices with positive determinant).
Hence, we can identify $\wt{\LG_A(\R)}$ with the quotient of
$\wt{\Sp}(2n,\R)$ by a connected subgroup, so $\wt{\LG_A(\R)}$ is simply connected.
It follows that in this case $\wt{\LG_A(\R)}$ is the
universal covering of the Lagrangian Grassmannian $\LG_A(\R)$.
\end{ex}


\subsection{Phase function}\label{phase-sec}

Since $\De:\BU(\R)\to\OO^*(D_A)$ is a $1$-cocycle (see Lemma \ref{De-coc-lem}), 
it defines a natural action of
the group $\BU(\R)$ (by holomorphic automorphisms) on the trivial $\C^*$-bundle over $D_A$.
We have constructed the central extension $U^\De\to \BU(\R)$ by $\Z$ in such a way
that $\De$ lifts to a $1$-cocycle of $U^\De$ with coefficients in $\OO(D_A)$. In other words,
we obtain the action of $U^\De$ on $D_A\times\C$ (respecting the structure of a $\C$-space), 
which we view as a universal covering of 
$D_A\times\C^*$ (in Sec.\ \ref{surface-sec} 
we will relate this covering to the Bridgeland's stability space in the case $\dim A=2$). 
Explicitly, this action is given by
\begin{equation}\label{DA-C-action}
(g,f)\cdot (\om,z)=(g(\om),z-f(\om)),
\end{equation}
where $(g,f)\in U^\De$ and $(\om,z)\in D_A\times\C$.

On the other hand, we have a transitive action of $U^\De$ on the $\Z$-covering 
$\wt{\LG(\R)}$ of $\LG(\R)$. 
By definition of this $\Z$-covering, we have a continuous function
$$\Bf_0: D_A\times\wt{\LG(\R)}\to \R: (\om, (L,f_L))\mapsto \Re f_L(\om).$$
We can extend it to a continuous function on $(D_A\times\C)\times\wt{\LG(\R)}$
setting
$$\Bf((\om,z),\wt{L})=\Re(z)+\Bf_0(\om,\wt{L}),$$
where $\wt{L}\in\wt{\LG(\R)}$.

\begin{lem}\label{equiv-lem} 
The function $\Bf$ is $U^\De$-invariant, i.e., 
for $\wt{g}\in U^\De$ and $(\si,\wt{L})\in (D_A\times\C)\times\wt{\LG(\R)}$
one has
$$\Bf(\wt{g}(\si),\wt{g}(\wt{L}))=\Bf(\si,\wt{L}).$$
\end{lem}

The proof is straightforward.


Now using the map $F\mapsto \wt{L}_F$ of Proposition \ref{Z-cov-prop}, 
we define the {\it phase function}
\begin{equation}\label{phi-si-eq}
(D_A\times\C)\times \ov{\SH}^{LI}(A)/\N^*\to \R:(\si,\wt{L})\mapsto \phi_\si(F):=\Bf(\si,\wt{L}_F),
\end{equation}
where $\si\in D_A\times\C$.
Note that we have
\begin{equation}\label{phi-plus-z}
\phi_{(\om,z)}(F)=\phi_{(\om,0)}(F)+\Re(z).
\end{equation}

In Sec.\ \ref{surface-sec} we will show that in the surface case the function $\phi_\si$ gives 
the phases of LI-objects with respect to
the Bridgeland's stability condition on $D^b(A)$ associated with $\si\in D_A\times\C$.
In the following theorem we check some of the properties of $\phi_\si$ that conform with the
conjecture that the corresponding stability condition exists in the higher-dimensional case as well.

\begin{thm}\label{phase-thm} The phase function $\phi_\si(F)$ satisfies the following properties.

\noindent
(i) This function is $\wt{\BU(\Q)}$-invariant, i.e.,
$$\phi_{\wt{g}(\si)}(\wt{g}(F))=\phi_\si(F),$$
where the action of $\wt{\BU(\Q)}$ on $D_A\times\C$ is induced by \eqref{DA-C-action}
via the homomorphism
$\iota:\wt{\BU(\Q)}\to U^\De$. In particular, for $n\in\Z$,
$$\phi_{\si}(F[n])=\phi_{\si}(F)+n.$$ 

\noindent
(ii) For $\si=(\om,z)$ and $F\in\ov{\SH}^{LI}(A)$ one has
\begin{equation}\label{phase-property}
\exp(\pi i z)\cdot \chi(\ell(\om),[F])\in \R_{>0}\cdot \exp(\pi i \phi_{\si}(F)),
\end{equation}
where $[F]\in \NN(A)\ot\R$ is the numerical class of $F$.
 
\noindent
(iii) For a semihomogeneous vector bundle $V_\phi$ associated with $\phi\in\NS(A)\ot\Q$ one
has
$$\phi_{(\om,z)}(V_\phi)=\Re(z)+\frac{1}{2\pi}\Arg(\deg(\om-\phi)),$$
where the branch of $\Arg(\deg(\cdot))$ is normalized by
$\Arg(\deg(iH))=-g\pi$ for ample $H$.

\noindent
(iv) for a pair of LI-objects $F_1$ and $F_2$ such that the corresponding
Lagrangians $L_{F_1}$ and $L_{F_2}$ in $A\times\hat{A}$ are transversal one has
\begin{equation}\label{main-phi-ineq}
\phi_{\si}(F_1)\le\phi_{\si}(F_2)+i(F_1,F_2),
\end{equation}
where $i(F_1,F_2)$ is the index of the pair $(F_1,F_2)$, i.e., 
the number such that $\Ext^i(F_1,F_2)=0$ for $i\neq i(F_1,F_2)$ (it exists by 
\cite[Cor.\ 3.2.12]{P-LIF}).
\end{thm}

\Pf . (i) The invariance follows from Lemma \ref{equiv-lem}. The second assertion follows
from this:
$$\phi_{(\om,z)}(F)=\phi_{(1,n)\cdot(\om,z)}(F[n])=\phi_{(\om,z-n)}(F[n])=\phi_{(\om,z)}(F[n])-n,$$
where in the last equality we used \eqref{phi-plus-z}.

\noindent
(ii) By part (i), 
the right-hand side of \eqref{phase-property} 
is invariant under the diagonal action of $\wt{\BU(\Q)}$ on $(\si,F)$.
We claim that the same is true for the left-hand side (modulo $\R_{>0}$). Indeed,
by Corollary \ref{chi-ell-cor}, for $\wt{g}=(g,f_g)\in U^\De$ we have 
$$\chi(\ell(\om),[F])=\exp(-\pi i f_g(\om))\cdot \chi(\ell(g(\om)),[\wt{g}(F)])\mod\R_{>0}$$
(recall that the map $F\mapsto [F]\mod\N^*$ is compatible 
with the projection $U^\De\to \Spin(\R)$
sending $(g,f_g)\in U^\De$ to $(g,\exp(-\pi i f_g))$, see \eqref{U-De-Spin-isom}). 
This immediately implies
that the left-hand side of \eqref{phase-property} is invariant modulo $\R_{>0}$ with
respect to the diagonal action of $\wt{\BU(\Q)}$ on
the pair $(\si,F)\in (D_A\times\C)\times\ov{\SH}^{LI}(A)/\N^*$.
 
Thus, it is enough to check the equality for $F=\OO_x$. We have
$\chi(\ell(\om),[\OO_x])=1$ for all $\om$. On the other hand, by definition of the map of
Proposition \ref{Z-cov-prop}, $\Bf_0(\om,\wt{L}_{\OO_x})=0$, so
$\phi_{(\om,z)}(\OO_x)=\Re(z)$.

\noindent
(iii) This follows from Proposition \ref{Z-cov-prop}(ii).

\noindent
(iv) By $\wt{\BU(\Q)}$-invariance of both parts with respect to the diagonal action on the pair
$(F_1,F_2)$, it is enough to consider the case when $F_2=\OO_x$. Note that in this case
the transversality assumption implies that $L_{F_1}=\Ga(\phi)$ for $\phi\in\NS(A)_\Q$,
so $F_1=V_\phi[n]$ for some $n\in\Z$, where $V_\phi$ is the simple semihomogeneous
bundle associated with $\phi$. Since $\phi_{(\om,z)}(\OO_x)=\Re(z)$, by part (iii),
the required inequality is equivalent to
$$\Arg(\deg(\om-\phi))\le 0,$$
where $\Arg(\deg(\cdot))$ is normalized by $\Arg(\deg(iH))=-g\pi$.
But this follows immediately from Lemma \ref{Arg-ineq-lem}(ii).
\ed


\begin{rem} 
By Lemma \ref{Arg-ineq-lem}(i), 
for $F_1=\OO$, $F_2=V_\phi$ and $\si=(iH,0)$, where $H$ is an ample class, 
the inequality \eqref{main-phi-ineq}
can be replaced by a stronger one:
$$\phi_{iH,0}(\OO)<\phi_{iH,0}(V_\phi)+\frac{i(\phi)}{2}.$$
However, this inequality is not invariant with respect to the group action considered above,
so it cannot be extended to the case of arbitrary $\si\in D_A\times\C$.
\end{rem}

The following property is also motivated by the picture with the stability conditions for $\dim A=2$
(see Sec.\ \ref{surface-sec} below).

\begin{prop}\label{Z-fibers-prop} 
The fibers of the map
$$\ZZ:D_A\times\C\to\Hom(\NN(A),\C):(\om,z)\mapsto \exp(\pi i z)\chi(\ell(\om),\cdot)$$
are exactly the orbits of the action of $2\Z\sub\C$ by translations on the second factor.
\end{prop}

\Pf . Suppose 
$$\exp(\pi i z)\chi(\ell(\om),[F])=\exp(\pi i z')\chi(\ell(\om'),[F])$$
for all $F$. Since $\chi(\ell(\cdot),[\OO_x])=1$, this implies
that $\exp(\pi i z)=\exp(\pi i z')$. Using the action of $2\Z$ we can assume that $z=z'$.
Now the fact that $\om=\om'$ follows from Corollary \ref{chi-cor}.
\ed

\begin{ex}\label{standard-stability-ex}
Recall that the standard stability condition on an elliptic curve has 
$Z(F)=-\deg(F)+i\rk(F)$ and semistable objects that are shifts of semistable bundles
and torsion sheaves. The corresponding phase function $\phi^{st}$ satisfies
$$\phi^{st}(F)=\phi_{(i,0)}(F)+1$$ 
for any semistable $F$. Indeed, this follows from the formulas
$$\phi^{st}(\OO_x)=1, \ \ \phi^{st}(V_{d/r})=\frac{\Arg(i-d/r)}{\pi},$$ 
where $V_{d/r}$ is the simple bundle of degree $d$ and rank $r$ and we normalize the
argument in the upper half-plane by $\Arg(i)=1/2$.
\end{ex}

\subsection{Stability conditions on abelian surfaces}\label{surface-sec}

In this section, assuming that $\dim A=2$
we will identify the action of $\iota(\wt{\BU(\Q)})\sub U^\De$ on
$D_A\times \C$
with the natural action of $\wt{\BU(\Q)}$ on the component $\Stab^{\dagger}$ of Bridgeland's
stability space $\Stab(A)$ of $D^b(A)$ described in \cite[Sec.\ 15]{Bridge-K3}. 

Recall that the stability space $\Stab(A)$ carries a natural continuous action of the group
$\wt{\GL}^+(2,\R)$, the universal cover of $\GL^+(2,\R)$, that can be described as the set
of pairs $(T,f)$, where $T\in\GL^+(2,\R)$ and $f:\R\to R$ is an increasing map with
$f(t+1)=f(t)+1$ such that the map induced by $T$ on $\R^2\setminus\{0\}/\R_{>0}\simeq \R/2\Z$ coincides with $f\mod 2\Z$. We use the left action of $\wt{\GL}^+(2,\R)$ on $\Stab(A)$: a pair $(T,f)$
maps a stability condition $(Z,\PP)$
to the stability $(T\circ Z, \PP')$, where $\PP'(t)=\PP(f^{-1}(t))$.
Note that $n\mapsto ((-1)^n, t\mapsto t+n)$ gives an embedding $\Z\to\wt{\GL}^+(2,\R)$ such that
$2\Z$ is the kernel of the projection to $\GL^+(2,\R)$.

Recall that for each $\om=i\a+\b\in D_A$ Bridgeland defined a stability condition on
$D^b(A)$ with the central charge
$$Z_\om(F)=-\chi(\ell(\om),[F])$$
and with each $\OO_x$ stable of phase $1$.
This defines a submanifold $V(A)\sub\Stab(A)$, isomorphic to $D_A$, which is a section of the action of 
$\wt{\GL}^+(2,\R)$ on a connected component $\Stab^{\dagger}(A)\sub\Stab(A)$, 
so that we have an isomorphism
\begin{equation}\label{Stab-answer}
V(A)\times\wt{\GL}^+(2,\R)\simeq\Stab^{\dagger}(A)
\end{equation}
(see \cite[Sec.\ 11, 15]{Bridge-K3}). Conjecturally, $\Stab^{\dagger}(A)=\Stab(A)$. Below we will
show that $\Stab^\dagger(A)$ contains all full stabilities (see Proposition \ref{ab-sur-stab-prop}). 

We have a natural embedding $\C^*=\GL(1,\C)\to\GL^+(2,\R)$ and the corresponding homomorphism
of universal coverings $\C\hra\wt{\GL}^+(2,\R)$ (where we use the map $\C\to\C^*:z\mapsto
\exp(\pi i z)$).
Hence, from the isomorphism \eqref{Stab-answer} we obtain an embedding
\begin{equation}\label{domain-stab-emb}
D_A\times\C\simeq V(A)\times\C\hra \Stab^{\dagger}(A).
\end{equation}
Note that the central charge corresponding to a point $(\om,z)\in D_A\times\C$ is
$$Z_{(\om,z)}(F)=-\exp(\pi i z)\chi(\ell(\om),[F]),$$
and the phase of $\OO_x$ with respect to this stability is 
\begin{equation}\label{Bridge-phase-O-x-eq}
\phi^{Br}_{(\om,z)}(\OO_x)=1+\Re(z)
\end{equation}
Recall that the non-empty fibers of the projection 
$$\ZZ:\Stab^{\dagger}(A)\to\Hom(\NN(A),\C)$$ 
are exactly the orbits of the action of $2\Z\sub\C\sub\wt{\GL}^+(2,\R)$
(see \cite[Thm.\ 15.2]{Bridge-K3}).
Hence, we have
\begin{equation}\label{D-C-preimage}
V(A)\times\C=\ZZ^{-1}(\ZZ(V(A)\times\C)).
\end{equation}
The image $\ZZ(V(A)\times\C)$ coincides with $\C^*\cdot \ell(D_A)\sub\NN(A)\ot\C$, where we
identify $\NN(A)\ot\C$ with $\Hom(\NN(A),\C)$ using $\chi(\cdot,\cdot)$.

Recall that we have an action of the group $\wt{\BU(\Q)}$ on $\Stab^{\bH}(A)$ defined using functors
$\Phi_g$ (see Def.\ \ref{stab-action-defi}). Also, note that by \cite[Cor.\ 3.5.2]{P-const}, we have an
inclusion $\Stab^{\dagger}(A)\sub\Stab^{\bH}(A)$
since all stabilities in $\Stab^{\dagger}(A)$ are full.

\begin{prop}\label{stab-action-prop}
The subset $V(A)\times\C\sub \Stab^{\bH}(A)$ is invariant with respect to the action of
$\wt{\BU(\Q)}$ and the induced action of $\wt{\BU(\Q)}$ on $V(A)\times\C\simeq D_A\times\C$
is exactly \eqref{DA-C-action}.
\end{prop}

\Pf . First, let us look at the action on central charges. 
Applying Corollary \ref{chi-ell-cor} to the element $\wt{g}=(g,\exp(-\pi i f))\in\Spin(\R)$
coming from an element $(g,f)=\iota(g')\in U^\De$ where $g'\in\wt{\BU(\Q)}$, we get
$$Z_{(\om,z)}(\hat{\rho}(g')^{-1}F)=Z_{(\om,z)}(\hat{\rho}(\wt{g})^{-1}[F])=Z_{\iota(g')\cdot(\om,z)}(F).$$
(see \eqref{DA-C-action}). In particular, the transformed central charge is still in $\C^*\cdot\ell(D_A)$.
Recall that the connected component $\Stab^{\dagger}(A)$ is characterized by the condition
that the central charge is in the $\GL^+(2,\R)$-orbit of $\ell(D_A)$ and $\OO_x$ are stable of the same
phase for all $x\in A$. Furthermore, by \cite[Lem.\ 12.2]{Bridge-K3}, is is enough to require
all $\OO_x$ to be semistable of the same phase (due to the absence of spherical objects---see
\cite[Lem.\ 15.1]{Bridge-K3}). In our case the condition on the central charge is satisfied by the 
above computation, and the semistability of $\OO_x$ follows from Proposition \ref{LI-semistable-prop}, 
so we get the inclusion $g'(V(A)\times\C)\sub\Stab^{\dagger}$.
Taking into account \eqref{D-C-preimage} we derive the required inclusion 
$$g'(V(A)\times \C)\sub V(A)\times\C\sub\Stab^{\dagger}.$$
Furthermore, we obtain that the action of $g'\in \wt{\BU(\Q)}$ on $D_A\times\C$ differs from
the action \eqref{DA-C-action} by the translation by an element in $2\Z\sub\C$.
Thus, the difference between the two action is given by a homomorphism
$\wt{\BU(\Q)}\to 2\Z$. Note that the element $1\in\Z\sub\wt{\BU(\Q)}$ acts on a stability
in $\Stab(A)$ by changing the central charge $Z$ to $-Z$ and adding $-1$ to all the phases.
Since this matches with its action on $D_A\times\C$ given by \eqref{DA-C-action}, 
the above homomorphism factors through a homomorphism $\BU(\Q)\to 2\Z$.
Next, we observe that 
the action of $\BP^-(\Q)\sub\wt{\BU(\Q)}$ preserves the phase of $\OO_x$ (see the proof of
Prop.\ \ref{Z-cov-prop}(i)). On the other hand, $\iota(\BP^-(\Q))\sub U^\De$ consists of elements
$(g,f)$ with $\Re(f)=0$ (see Lemma \ref{U-De-splitting-lem}), 
so taking into account the formula \eqref{Bridge-phase-O-x-eq} we deduce that
the homomorphism $\BU(\Q)\to 2\Z$ is trivial on $\BP^-(\Q)$. It remains to apply Lemma 
\ref{P-hom-lem}.
\ed

\begin{cor}\label{stab-cor} 
There is a transitive continuous action of $U^\De\times \wt{\GL}_2^+(\R)$ on
$\Stab^{\dagger}(A)$, extending the action of $\wt{\BU(\Z)}$ (coming from autoequivalences of
$D^b(A)$) and the standard action of
$\wt{\GL}_2^+(\R)$.
\end{cor}

\Pf . This follows from the identification \eqref{Stab-answer} and from the transitivity
of the action of $U^\De$ on $D_A\times\C$. Note that our action
 of $\wt{\BU(\Q)}$ on $\Stab^{\bH}(A)$ extends the standard action of $\wt{\BU(\Z)}$
by autoequivalences of $D^b(A)$.
\ed

\begin{thm}\label{surface-thm} 
For any $\si=(\om,z)\in D_A\times\C$ and any LI-object $F\in D^b(A)$, let
$\phi_{\si}^{Br}(F)$ be the phase of $F$ with respect to the corresponding Bridgeland's stability condition.
Then 
$$\phi_{\si}^{Br}(F)=\phi_{\si}(F)+1,$$ 
where the function $\phi_{\si}$ is given by \eqref{phi-si-eq}.
\end{thm}

\Pf . The assertion is true for $F=\OO_x$. Also Theorem \ref{phase-thm}(i) together with
Proposition \ref{stab-action-prop} imply that both
sides are invariant with respect to the action of $\wt{\BU(\Q)}$ on the pair $(\si,F)$.
It remains to use transitivity of the action of
$\wt{\BU(\Q)}$ on $\ov{\SH}^{LI}(A)/\N^*$.
\ed

We finish our consideration of abelian surfaces by observing that Proposition \ref{LI-semistable-prop}
implies the following description of the component $\Stab^{\dagger}(A)$. 

\begin{prop}\label{ab-sur-stab-prop} For an abelian surface $A$ the component $\Stab^\dagger(A)$
of the stability space consists of all full stabilities on $D^b(A)$.
\end{prop}

\Pf . Let $\sigma$ be a full stability on $D^b(A)$.
By Proposition \ref{LI-semistable-prop}, any skyscraper sheaf $\OO_x$ is $\si$-semistable.
Using Lemmas 12.2 and 15.1 of \cite{Bridge-K3} we deduce that $\OO_x$ is in fact $\si$-stable.
Now the proof of \cite[Prop.\ 3.1.4]{Bridge-K3} shows that any full stability belongs to $\Stab^\dagger(A)$
(one has to run this proof with $U(X)$ being the set of all full stabilities, which is an open subset
of $\Stab(A)$).
\ed

\subsection{Mirror symmetry and phases}\label{mirror-sec}

In the case when $A=E^n$, where $E$ is an elliptic curve without complex mutliplication,
we can interpret the phase function of Sec.\ \ref{phase-sec} in terms of the Fukaya category
of the mirror dual abelian variety.

Let $E=\C/(\Z+\tau\Z)$ be an elliptic curve over $\C$ and let
$\La$ be a free $\Z$-module of rank $n$. 
We set 
$$A=\La\ot\C/(\La\ot (\Z+\tau\Z))\simeq E^n,$$
so that we have a natural isomorphism
$$\Ga_A:=H_1(A,\Z)\simeq \La\oplus\La,$$
where the second summand corresponds to $\La\ot \tau$.
The natural polarization of $E$ given by the hermitian form 
$H_\tau(z_1,z_2)=\frac{z_1\ov{z_2}}{\Im\tau}$ induces an isomorphism
$$\hat{A}\simeq \La^*\ot\C/(\La^*\ot (\Z+\tau\Z),$$
where $\La^*=\Hom_\Z(\La,\Z)$.

Assuming that $E$ has no complex multiplication we obtain identifications
$\End(A)\simeq\End_\Z(\La)$, $\Hom(A,\hat{A})\simeq\Hom_\Z(\La,\La^*)$, and
$\NS(A)\simeq \Hom_\Z(\La,\La^*)^+$ (the latter group consists of symmetric homomorphisms).
Thus, for a field $F\supset\Q$ 
we can identify $\NS(A)\ot F$ with the space of symmetric bilinear forms on 
$\La\ot F$. The ample cone in $\NS(A)\ot\R$ consists of positive-definite forms. Thus,
$D_A\sub\NS(A)\ot\R$ is the Siegel's half-space consisting of symmetric bilinear forms on
$\La\ot\C$ with positive-definite imaginary part.

According to \cite{GLO} (see also \cite[Sec.\ 6.5]{P-ab-var}), one can view the 
abelian variety $B$ associated with
an element $\om=\om_A\in D_A\simeq\HP_n$ as a mirror dual to $(A,\om_A)$. More precisely, let us set 
$$\Ga_B=\La^*\oplus\La, \ \ B=\Ga_B\ot \R/\Ga_B,$$
and define the complex structure on $\Ga_B\ot \R$ via the isomorphism
\begin{equation}\label{kappa-om-eq}
\kappa_\om:\Ga_B\ot \R\to \La^*\ot\C: (\la^*,\la)\mapsto \la^*-\om(\la),
\end{equation}
where we view $\om$ as an element of $\Hom(\La,\La^*\ot\C)^+$.
Note that there is an isomorphism $B\simeq\La^*\ot\C/(\La^*+\om\La)$ (however,
the corresponding identification of $H_1(B,\Z)$ with $\Ga_B$ differs from the original one 
by the sign on the summand $\La\sub\Ga_B$).
We have a natural principal polarization $\phi_0:B\wt{\to}\hat{B}$ given on homology lattices by
\begin{equation}\label{principal-polarization-B-eq}
\Ga_B\to\Ga_B^*: (\la_0^*,\la_0)\mapsto ((\la^*,\la)\mapsto \la^*(\la_0)-\la_0^*(\la).
\end{equation}
Similarly, the natural isomorphism $\La^*\oplus\La^*\simeq\Ga_{\hat{A}}\simeq\Ga_A^*$ corresponds
to the pairing
$$(\La^*\oplus\La^*)\times\Ga_A\to\Z: ((\la_1^*,\la_2^*),(\la_1,\la_2))\mapsto \la_1^*(\la_2)-\la_2^*(\la_1).$$

Let us define an isomorphism of orthogonal lattices
$$\ga:\Ga_A\oplus\Ga_{\hat{A}}\to \Ga_B\oplus\Ga_B\simeq\Ga_B\oplus\Ga_{\hat{B}}:
(\la_1,\la_2,\la_1^*,\la_2^*)\mapsto (\la_2^*,\la_2,\la_1^*,\la_1).$$

\begin{prop}\label{mirror-prop} The isomorphism $\ga$ induces a mirror duality in the sense of
\cite[Sec.\ 9]{GLO} between the pairs $(A,\om_A)$ and $(B,\om_B)$ for $\om_B=\tau\cdot\phi_0$,
where $\phi_0\in\Hom(B,\hat{B})^+$ is the principal polarisation defined above.
\end{prop}

\Pf . By definition, we have to check that the operator of complex structure on 
$(\Ga_B\oplus\Ga_B)\ot\R$ corresponds under $\ga$ to 
$$I_{\om_A}=\left(\begin{matrix} \a^{-1}\b & -\a^{-1}\\ \a+\b\a^{-1}\b & -\b\a^{-1}\end{matrix}\right)\in
\BU_A(\R),$$
where $\om_A=i\a+\b$, and we view $\BU_A(\R)$ as a subgroup in automorphisms of
$(\Ga_A\oplus\Ga_{\hat{A}})\ot\R$, and similarly, that the complex structure on 
$(\Ga_A\oplus\Ga_{\hat{A}})\ot\R$ corresponds to $I_{\om_B}$.
Both facts are checked by a straightforward computation (cf.\ \cite[Prop.\ 9.6.1]{GLO}).
\ed

Recall that the variety $\LG_A=\LG_{E^n}$ is naturally identified with the Lagrangian Grassmannian associated
with the symplectic lattice $\La^*\oplus\La$. Thus, a Lagrangian subvariety $L\sub A\times\hat{A}$,
viewed as a point of $\LG(\Q)$, corresponds to a Lagrangian $\Z$-submodule 
$\Pi(L)\sub \La^*\oplus\La=\Ga_B$,
so that $\Ga_L=H_1(L,\Z)\simeq\Pi(L)\oplus\Pi(L)\sub\Ga_A\oplus\Ga_{\hat{A}}$. 
Hence, from \eqref{de-L-Lag-eq} we get
\begin{equation}\label{de-kappa-om-eq}
\de(L)(\om)=\det(\kappa_\om|_{\Pi(L)})^2\mod\R_{>0},
\end{equation}
where we view $\kappa_\om|_{\Pi(L)}$ as an element in $\Hom(\Pi(L),\La^*)\ot\C$
and define $\det^2$ using some bases in $\Pi(L)$ and $\La^*$.

Similarly, a point $L$ of $\LG_A(\R)$ corresponds to a real Lagrangian subspace 
$\Pi_\R(L)\sub\Ga_B\ot\R$ and
the formula \eqref{de-kappa-om-eq} still holds (with $\Pi_\R(L)$ instead of $\Pi(L)$).
Recall that we have a double covering $\LG^{\spin}(A,\R)\to\LG(\R)$ consisting of pairs
$(L,\varphi)\in \LG_A(\R)\times\OO(D_A)/\R_{>0}$ such that $\varphi^2=\de(L)$,
so that the group $\Spin(\R)$ acts on $\LG^{\spin}(A,\R)$ (see Def.\ \ref{LG-spin-defi}). 
In our case there is a splitting $\BU(\R)\to\Spin(\R)$ (see Remark \ref{spin-rem}.1), 
so we have an action of $\BU(\R)$ on $\LG^{\spin}(A,\R)$ given by
$$g\cdot (L,\varphi)=(gL,\varphi'), \ \text{where }\ \varphi'(g(\om))=\varphi(\om)\cdot \det(a+b\om)^{-1}.
$$

We claim that $\LG^{\spin}(A,\R)\to\LG_A(\R)$ is in fact the 
natural double covering corresponding to a choice of orientation on a Lagrangian subspace in 
$\Ga_B\ot\R$. Indeed, let us fix an orientation $\eps\in\bigwedge^n(\La)$.
Then a choice of a square root $\varphi=\sqrt{\de(L)}\in\OO^*(D_A)/\R_{>0}$ for $L\in\LG_A(\R)$
induces an orientation on $\Pi_\R(L)\sub\Ga_B\ot\R$ 
as follows. By formula \eqref{de-kappa-om-eq}, for each $\om$ the non-zero element
$$\varphi(\om)^{-1}\cdot\det(\kappa_\om|_{\Pi_\R(L)})\in
{\bigwedge}^n(\La^*)\ot{\bigwedge}^n(\Pi_\R(L))^{-1}\ot_\R \C,$$
depending continuously on $\om$,
belongs to the $\R$-subspace $\bigwedge^n(\La^*)\ot\bigwedge^n(\Pi_\R(L))^{-1}$. Thus, we get an isomorphism
$${\bigwedge}^n(\Pi_\R(L))\simeq{\bigwedge}^n(\La)^{-1}\ot\R$$
and we define the orientation $\mu_{\varphi,\eps}\in \bigwedge^n(\Pi_\R(L))$ so that it corresponds to $\eps^{-1}$ under this isomorphism, i.e.,
$$\varphi(\om)^{-1}\cdot\det(\kappa_\om|_{\Pi_\R(L)})\cdot \mu_{\varphi,\eps}=\eps^{-1}.$$

Let us associate with $L\in\LG_A(\Q)$ the real subtorus in $B$ by setting
$$T_L=\Pi(L)\ot\R/\Pi(L)\sub \Ga_B\ot \R/\Ga_B=B.$$
Note that $T_L$ is Lagrangian with respect to the translation-invariant symplectic structure on
$B$ corresponding to the standard symplectic structure on $\Ga_B=\La^*\oplus\La$
(i.e., this symplectic structure on $B$ comes from the principal polarization $\phi_0$).
As we have shown above, 
a lifting of $L$ to a point $(L,\varphi)\in\LG^{\spin}(A,Q)$ gives rise to an orientation of $T_L$.

Since $\LG^{\spin}(A,\Q)=\wt{\LG_A(\Q)}/2\Z$, the map \eqref{Z-cov-bijection} induces a map
$$\ov{\SH}^{LI}\to\LG^{\spin}(A,\Q).$$
By Lemma \ref{NS-section-lem}, the composition
$$\NS(A)\ot\Q\to\ov{\SH}^{LI}\to\LG^{\spin}(A,\Q): \phi\mapsto V_\phi\mapsto (\Ga(\phi),\varphi)$$
corresponds to the choice of the square root $\varphi(\om)=\det(\phi-\om)$,
where we use dual bases of $\La$ and $\La^*$ to compute the determinant
(see also Ex.\ \ref{ell-power-ex}). 
The corresponding orientation on $T_{\Ga(\phi)}$ is induced by
the isomorphism $\Pi(\Ga(\phi))\ot\R \simeq\La\ot\R$ and the orientation $\eps$ of $\La\ot\R$.

Let $\Om_{\om,\eps}$ denote the holomorphic volume form on $B$ defined by
$$\Om_{\om,\eps}=\kappa_\om^*(\eps),$$
where we view $\eps\in\bigwedge^n(\La)\sub\bigwedge^n(\La)\ot\C$ as an $n$-form on 
$\La^*\ot\C$ and use an isomorphism \eqref{kappa-om-eq}.

\begin{thm}\label{integral-thm} 
For an endosimple LI-object $F\in D^b(A)$ one has
\begin{equation}\label{integral-eq}
\chi(\ell(\om),[F])=\int_{[T_L]}\Om_{\om,\eps}
\end{equation}
where $L=L_F$, and $T_L$ is equipped with the orientation $\mu_{\varphi_F,\eps}$
coming from the point $(L_F,\varphi_F)\in\LG^{\spin}(A,\Q)$ associated with $F$.
\end{thm}

\Pf . Note that shifting $F$ by $[1]$ changes the orientation of $T_L$, so the assertions for
$F$ and $F[n]$ are equivalent.

First, let us prove \eqref{integral-eq} in the case when $L_F$ is transversal to $\{0\}\times\hat{A}$,
i.e., when $F=V_\phi$ is the semihomogeneous bundle corresponding to
$\phi\in\Hom(A,\hat{A})^+\ot\Q\simeq\Hom(\La,\La^*)^+\ot\Q$ (and $L_F=\Ga(\phi)\sub A\times\hat{A}$).
Recall that $\rk V_\phi=\deg(L_F\to A)^{1/2}$ (see \eqref{rk-V-phi-eq}).
For $K=\Q$ or $\R$ let $\Ga_{K}(\phi)\sub(\La^*\oplus\La)\ot K$ be the graph of $\phi$ viewed as a map of $K$-vector spaces (i.e., $\Ga_{K}(\phi)=H_1(L_F,K)$) and set
$$\Ga_{\Z}(\phi):=\Ga_{\Q}(\phi)\cap(\La^*\oplus\La),$$
so that $T_L=\Ga_{\R}(\phi)/\Ga_\Z(\phi)$. We also denote by $i_\phi:\Ga_\Z(\phi)\to\La^*\oplus\La$
the natural embedding. The orientation on $T_L$ is induced by the natural isomorphism
$\Ga_{\R}(\phi)\simeq\La\ot\R$ and by the orientation of $\La\ot\R$ given by $\eps$.
The cycle $[T_L]$ in $H_n(A)\simeq \La^n(\La^*\oplus\La)$ is the image
of the positive generator $\mu\in\bigwedge^n(\Ga_\Z(\phi))$ under the map 
$${\bigwedge}^n(i_\phi):{\bigwedge}^n(\Ga_\Z(\phi))\to{\bigwedge}^n(\La^*\oplus\La).$$
Note also that the integration map
$$H_n(A)\to\C:\ga\mapsto \int_\ga \Om_{\om,\eps}$$
is identified with 
$${\bigwedge}^n(\kappa_\om):{\bigwedge}^n(\La^*\oplus\La)\to{\bigwedge}^n(\La^*)\ot\C\simeq\C,$$
where the last isomorphism is given by $\eps$.
Hence, $\int_{[T_L]}\Om_{\om,\eps}=\de(\mu)\cdot\eps$, where
$\de\in\bigwedge^n(\La^*)\ot\bigwedge^n(\Ga_\Z(\phi))^{-1}\ot\C$ is the determinant of the composition
$$\Ga_\Z(\phi)\rTo{i_\phi}\La^*\oplus\La\rTo{\kappa_\om}\La^*.$$
The projection $p_2:\Ga_\Z(\phi)\to \La$ is an embedding of index $\deg(L_F\to A)^{1/2}$, so the commutative diagram
\begin{diagram} 
\Ga_\Z(\phi)&\rTo{\kappa_\om i_\phi}&\La^*\ot\C\\
\dTo{p_2}&&\dTo{\id}\\
\La&\rTo{\phi-\om}&\La^*\ot\C
\end{diagram}
implies that
$$\de(\mu)\cdot\eps=\det(\phi-\om)\cdot\deg(L_F\to A)^{1/2}=\chi(\ell(\om),\ell(\phi))\cdot\rk(F)=
\chi(\ell(\om),[F]),
$$ 
where the last equality follows from Lemma \ref{ell-lem}.

Next, we will check that \eqref{integral-eq} is compatible with the action of the group
$\BU(\Z)$ on $[F]$, $\om$ and on $B$, where we use the natural symplectic action of $\BU(\Z)$
on $\Ga_B=\La\oplus\La^*$ and the splitting $\BU(\Z)\to\BU(\Z)^{\spin}$ of the spin-covering (see 
Remark \ref{spin-rem}.1). Namely, for
$g=\left(\begin{matrix} a & b \\ c & d\end{matrix}\right)\in \BU(\Z)$ the relation
$$
\left(\begin{matrix} -\om & \id_{\La^*}\end{matrix}\right)\cdot g^{-1}=
(a+b\om)^*\cdot \left(\begin{matrix} -g(\om) & \id_{\La^*}\end{matrix}\right)
$$
leads to a commutative diagram
\begin{equation}\label{kappa-diagram}
\begin{diagram}
\Ga_B &\rTo{\kappa_\om}&\La^*\ot \C\\
\dTo{g}&&\uTo{}{(a+b\om)^*}\\
\Ga_B &\rTo{\kappa_{g(\om)}}&\La^*\ot\C
\end{diagram}
\end{equation}
Hence, we have 
$$g^*\Om_{g(\om),\eps}=\det(a+b\om)^{-1}\cdot \Om_{\om,\eps},
$$
which implies that
$$\int_{g[T_L]}\Om_{g(\om),\eps}=\det(a+b\om)^{-1}\cdot \int_{[T_L]}\Om_{\om,\eps}.
$$
Also, the diagram \eqref{kappa-diagram} gives the equation
\begin{equation}\label{kappa-om-L-eq}
\kappa_\om|_{\Pi(L)}=(a+b\om)^*\circ \kappa_{g(\om)}|_{\Pi(gL)}\circ g|_{\Pi(L)}
\end{equation}
in $\Hom(\Pi(L),\La^*)\ot\C$.
Since $g\cdot (L,\varphi)=(gL,\varphi')$, where $\varphi'(g(\om))=\varphi(\om)\det(a+b\om)^{-1}$,
passing to determinants in \eqref{kappa-om-L-eq} we obtain that
the orientation $\mu_{\varphi',\eps}$ of $\Pi(gL)\ot\R$ corresponds to $\mu_{\varphi,\eps}$
under the isomorphism $\Pi(L)\to\Pi(gL)$ given by $g$. Hence, the class $g[T_L]$ is exactly
the fundamental class of $T_{gL}$ associated with the orientation coming from $g[F]$.
On the other hand, by Corollary \ref{chi-ell-cor},
$$\chi(\ell(\om),[F])=\chi(\ell(g(\om)),g[F]),$$
since for $g\in \BU(\Z)$ the operator $\hat{\rho}(g)$ is simply the map induced by any autoequivalence of $D^b(A)$ compatible with the canonical lifting of $g$ to $\BU(\Z)^{\spin}$.

Finally, applying Proposition \ref{Lag-auteq-prop} and using the $\BU(\Z)$-invariance, we see
that the general case of \eqref{integral-eq} follows from
the case when $L_F$ is transversal to $\{0\}\times\hat{A}$ considered above.
\ed

\begin{rem} Note that since $L_F$ is equipped with the lifting $\wt{L}_F$ to the universal
covering of the Lagrangian Grassmannian (see Ex.\ \ref{ell-power-ex}), the Lagrangian torus
$T_L$ has a structure of a {\it graded Lagrangian} (see \cite{Seidel}). The
corresponding choice of a phase of $\int_{T_L}$ obtained from Theorems \ref{phase-thm} and 
\ref{integral-thm}
comes from Kontsevich's description of a grading on a Lagrangian (see \cite[Ex.\ 2.9]{Seidel}).
\end{rem}

\section{Quasi-standard $t$-structures and Fourier-Mukai partners}\label{t-str-FM-sec}

\subsection{Quasi-standard $t$-structures}

The $\Z$-covering $\wt{\LG(\Q)}\to\LG(\Q)$ appears also naturally when considering $t$-structures.
Let $\TT(A)$ be the set of $\bH$-invariant $t$-structures on $D^b(A)$.
We identify $\TT(A)$ with the set of cores of such $t$-structures, so we view elements of $\TT(A)$
as abelian subcategories $\AA\sub D^b(A)$.

\begin{thm}\label{t-str-thm} 
(i) There is a natural $\wt{\BU(\Q)}$-equivariant embedding 
$$\wt{\LG(\Q)}\to \TT(A): \wt{L}\mapsto \AA_{\wt{L}},$$
which is uniquely characterized by the condition
$$\AA_{(0:\phi_0),0}=\Coh(A).$$
The LI-functor $\Phi_{\wt{g}}:D^b(A)\to D^b(A)$ 
corresponding to $\wt{g}\in\wt{\BU(\Q)}$ (defined up to $\bH$---see Sec.\ \ref{LI-obj-sec}) satisfies
$$\Phi_{\wt{g}}(\AA_{\wt{L}})\sub \AA_{\wt{g}\wt{L}}.$$

\noindent
(ii) For an LI-object $F$ and $\wt{L}\in\wt{\LG(\Q)}$
one has $F[-i(\wt{L}_F,\wt{L})]\in\AA_{\wt{L}}$, where $i(\cdot,\cdot)\in\Z$
is the unique $\wt{\BU(\Q)}$-equivariant function on $\wt{\LG(\Q)}\times \wt{\LG(\Q)}$
such that for $\phi_1,\phi_2\in\NS(A)\ot\Q$
one has
$$i(\wt{L}_{V_{\phi_1}},\wt{L}_{V_{\phi_2}})=i(\phi_2-\phi_1),$$
provided $\phi_2-\phi_1$ is nondegenerate (recall that $\wt{L}_{V_\phi}$ is given by 
\eqref{wt-L-V-phi-eq}).
\end{thm}

\Pf . (i) Recall that the action of $\wt{\BU(\Q)}$ on $\wt{\LG(\Q)}$ is transitive,
and the stabilizer subgroup of the point $((0:\phi_0),0)$ is $\BP^-(\Q)$, lifted to $\wt{\BU(\Q)}$
as described in Corollary \ref{lifting-cor}.
Thus, it suffices to check that $\BP^-(\Q)$ preserves the standard $t$-structure.
But this immediately follows from the description of the functors corresponding to elements of
$\BP^-(\Q)$ (see Prop.\ \ref{lifting-prop}).

\noindent
(ii) The fact that every 
LI-sheaf is cohomologically pure with respect to each $t$-structure constructed in (i) follows from
Theorem \ref{LI-sh-thm}. Uniqueness of the $\wt{\BU(\Q)}$-equivariant index function $i(\cdot,\cdot)$
follows from Proposition \ref{Lag-auteq-prop}. It remains to find the number $i=i(\phi_1,\phi_2)$ such that
$$V_{\phi_1}[-i]\in\AA_{\wt{L}_{V_{\phi_2}}}.$$
Let $g=\left(\begin{matrix} 1 & \phi_2^{-1} \\ 0 & 1\end{matrix}\right)$.
Then by formula \eqref{Phi-g-L-eq}, we have 
$$\Phi_g(\OO_x)=V_{\phi_2}\mod\N^*$$
(there is no shift in this case since the kernel $S(g)$ is a vector bundle).
It follows that 
$$\Phi_g(\Coh(A))\sub\AA_{\wt{L}_{V_{\phi_2}}}.$$
Note that $\Ga(\phi_1)=g(\Ga(\phi))$, where
$$\phi=\phi_1(1-\phi_2^{-1}\phi_1)^{-1}.$$
Hence, using \eqref{Phi-g-L-eq} and \eqref{la-g-phi-eq} we obtain
$$\Phi_g(V_\phi)=V_{\phi_1}[-i(\phi_2+\phi)] \mod \N^*,$$
so denoting $\phi_3=1-\phi_2^{-1}\phi_1$ we obtain
$$i=i(\phi_2+\phi_1\phi_3^{-1})=i(\phi_3(\phi_2\phi_3+\phi_1))=
i(\phi_3\phi_2)=i(\phi_2-\phi_1)$$
as claimed.
\ed

\begin{defi} We will refer to $t$-structures on $D^b(A)$ constructed in the above theorem
as {\it quasi-standard $t$-structures}.
\end{defi}

\begin{prop} Let $A$ and $B$ be abelian varieties, and let
$\eta: X_A\to X_B$ be a symplectic isomorphism in $\Ab_\Q$ (i.e., up to isogeny).
Then the map $\eta_*:\LG_A(\Q)\to \LG_B(\Q)$ induced by $\eta$ extends to a $\Z$-equivariant
map 
$$\wt{\eta}_*:\wt{\LG_A(\Q)}\to\wt{\LG_B(\Q)}$$
which is compatible with the quasi-standard $t$-structures, i.e., for every
$\wt{L}\in\wt{\LG_A(\Q)}$ 
the LI-functor $\Phi_\eta$ associated with $\eta$ (defined up to $\bH$)
satisfies 
\begin{equation}\label{t-exact-eq}
\Phi_\eta(\AA_{\wt{L}})\sub \AA_{\wt{\eta}_*\wt{L}}.
\end{equation}
\end{prop}

\Pf .
Note that $B$ is isogenous to $A$, i.e., there exists an isomorphism $f:A\to B$ in $\Ab_\Q$.
Let $\eta_0:X_A\to X_B$ be the induced symplectic isomorphism in $\Ab_\Q$. 
We also have natural compatible isomorphisms induced by $f$:
$$\BU_{X_A}\to \BU_{X_B}, \ \ \ \wt{\BU_{X_A}(\Q)}\to \wt{\BU_{X_B}(\Q)}$$
$$\eta_{0*}:\LG_A(\Q)\to\LG_B(\Q), \ \ \ \wt{\eta}_{0*}:\wt{\LG_A(\Q)}\to\wt{\LG_B(\Q)}.$$
Furthermore, it is easy to see that the $t$-exactness \eqref{t-exact-eq} holds for $\wt{\eta}_{0*}$
and the functor $\Phi_{\eta_0}$ which is the composition of the pull-back and the push-forward
under isogenies (this is proved similarly to Prop.\ \ref{lifting-prop}(ii)).
Now let $g_\eta\in \BU(\Q)$ be the unique element such that 
$$\eta=\eta_0\circ g_\eta.$$
Choose any element $\wt{g}_\eta\in\wt{\BU(\Q)}$ over $g_\eta$ and define
$$\wt{\eta}_*:\wt{\LG_A(\Q)}\to\wt{\LG_B(\Q)}: \wt{L}\mapsto \wt{\eta}_{0*}(\wt{g}_\eta(\wt{L})).$$
By Theorem \ref{t-str-thm}(i), the required assertion follows for the functor
$\Phi_{\eta_0}\circ\Phi_{\wt{g}_\eta}$. By \cite[Thm.\ 3.2.11]{P-LIF}, 
its  $\bH$-equivalence class 
differs from $\Phi_\eta[n]$ by an action of $\N^*$ (one has to use also \cite[Prop.\ 2.4.7(ii)]{P-LIF}
as in the proof of \cite[Thm.\ 3.3.4]{P-LIF}). 
Changing $\wt{\eta}_*$ using the action of $n\in\Z\sub\wt{\BU(\Q)}$
on $\wt{\LG_A(\Q)}$, we get the required compatibility \eqref{t-exact-eq}.
\ed

\begin{rems}
1. The quasi-standard $t$-structure associated with $\wt{L}_F\in\wt{\LG(\Q)}$ has a simple characterization in terms of the LI-object $F$ (defined up to $\bH$-equivalence).
Namely, the corresponding subcategory $D^{\le 0}\sub D^b(A)$ consists of all
$X\in D^b(A)$ such that $\Hom^i(X, T_{x,\xi}(F))=0$ for $i<0$ and all $(x,\xi)\in A\times\hat{A}$.
Indeed, using $\wt{\BU(\Q)}$-action this reduces to the characterization of the standard
subcategory $D^{\le 0}$ by the above condition, where $F$ is a nonzero torsion sheaf.

\noindent 2.
In the case of an elliptic curve all the quasi-standard $t$-structures are obtained from the
standard one by tilting (up to a shift). More precisely, let $P(\cdot)$ be the slicing associated
with the standard stability on $D^b(E)$ for an elliptic curve $E$, so that $P((0,1])=\Coh(E)$
(see Ex.\ \ref{standard-stability-ex}).
Then the quasi-standard $t$-structure associated with $\phi\in\NS(E,\Q)\simeq\Q$ (lifted to
$\wt{\LG_E(\Q)}$ by \eqref{NS-LG-section-eq}) is 
$P((\frac{\Arg(i-\phi)}{\pi}-1,\frac{\Arg(i-\phi)}{\pi}])$.
Note that this construction extends to irrational numbers $\phi\in\R$ and for $k=\C$ the corresponding
hearts are equivalent to the categories of holomorphic bundles on noncommutative $2$-tori
(see \cite{PS}, \cite{P-nc-tori}). We conjecture that this connection between
quasi-standard $t$-structures and  noncommutative tori
extends to the higher-dimensional case (the corresponding
equivalence of derived categories is established in \cite{Block}). 
Namely, to every point of $\wt{\LG_A(\R)}\setminus\wt{\LG_A(\Q)}$ there should correspond a $t$-structure on $D^b(A)$
(in a way compatible with the action of $\wt{\BU(\Q)}$) whose heart is equivalent
to the category of holomorphic bundles on the corresponding noncommutative torus.
\end{rems}


\subsection{Fourier-Mukai partners}\label{FM-sec}

Recall that the set of Fourier-Mukai partners ({\it FM-partners} for short)
of a smooth projective variety $X$ is
defined as
$$\FM(X)=\{Y \text{ smooth projective } |\ D^b(Y)\simeq D^b(X)\}/\text{isomorphism}.$$
For an abelian variety $A$ we can also define the subset
$\FM^{ab}(A)\sub\FM(A)$ by considering only FM-partners among abelian varieties.
In characteristic zero it is known that $\FM^{ab}(A)=\FM(A)$ (see
\cite{HN}).

Recall that if $B$ is a FM-partner of $A$ then any equivalence $D^b(A)\simeq D^b(B)$
is given by the LI-kernel associated with a Lagrangian correspondence $(L(\eta),\a)$
extending a symplectic isomorphism $\eta: X_A\simeq X_B$ 
(see \ref{LI-obj-sec}). 
The $\BU(\Z)$-orbit of the Lagrangian $(\eta_*)^{-1}(0\times \hat{B})\in \LG_A(\Q)$ does not
depend on a choice of an equivalence $D^b(A)\simeq D^b(B)$. 

\begin{prop}\label{FM-prop} 
The above construction gives an embedding 
\begin{equation}\label{FM-map}
\FM^{ab}(A)\hra\LG_A(\Q)/\BU(\Z).
\end{equation}
The image consists of orbits of Lagrangian subvarieties $L\sub X_A$ for which there exists
a Lagrangian subvariety $L'\sub X_A$ such that $L\cap L'=0$.
\end{prop}

\Pf . The first assertion is immediate since the Lagrangian subvariety 
$(\eta_*)^{-1}(0\times\hat{B})\sub X_A$ corresponding to $B$ is isomorphic to $\hat{B}$.
For the second we observe that if we have a Lagrangian $L'\sub X_A$ such that $L\cap L'=0$
then the get an isomorphism $L\times L'\to X_A$ and also $L'\simeq X_A/L\simeq \hat{L}$,
which leads to a symplectic isomorphism $L\times\hat{L}\simeq X_A$, so that
$B=\hat{L}$ is a FM-partner of $A$.
\ed

\begin{rem}\label{finite-rem}
The set $\LG_A(\Q)/\BU(\Z)=\BU(\Z)\backslash \BU(\Q)/\BP^-(\Q)$ is known to be finite
(see \cite[Thm.\ 6]{Go}). Note that this set is also in bijection with the set of endosimple
LI-objects in $D^b(A)$ up to the action of exact autoequivalences of $D^b(A)$ (as follows
from Prop.\ \ref{LG-map}).
\end{rem}

Here is an example of a situation when the embedding of Proposition \ref{FM-prop} is a
bijection.

\begin{prop}\label{class-prop} 
Assume that $A$ is principally polarized and $\End(A)=R$ is the ring of integers in
a totally real number field $F$
(so the Rosati involution on $F$ is trivial). Then the map \eqref{FM-map} is a bijection, and 
$$|\FM^{ab}(A)|=|\LG_A(\Q)/\BU(\Z)|=h_R,$$
where $h_R$ is the class number of $R$.
\end{prop}

\Pf . First, we observe that in this case the set $\LG_A(\Q)$ consists of all subvarieties in 
$X_A=A\times A$, isogenous to $A$. We claim that all such subvarieties $L\sub X_A$ are direct summands. Indeed, $L$ is an image of the morphism $A\to A^2$ associated with 
a pair $(a,b)\in R^2\setminus\{(0,0)\}$.
Consider the exact sequence
$$0\to I'\to R^2\to I\to 0,$$
where $I=(a,b)\sub R$. This sequence splits since $I$ is a projective $R$-module.
Hence, there is a corresponding split exact sequence of abelian varieties
$$0\to A^I\to A^2\to A^{I'}\to 0,$$
where we use the natural functor $M\to A^M$ from $R$-modules to commutative group schemes
with $A^M(S)=\Hom_R(M,A(S))$ (see \cite{Giraud})
Since $A^I$ is exactly the image of the map $(a,b):A\to A^2$, this proves our claim.

It remains to check that the orbits of $\SL_2(R)$ on the projective line $\P^1(F)$ are in bijection
with the ideal class group $\Cl(R)$. We have a well defined map 
$$\P^1(F)/\SL_2(R)\to\Cl(R)$$
sending $(a:b)$ with $a,b\in R$ to the class of the ideal $(a,b)$.
This map is surjective since every nonzero ideal in $R$ is generated by two elements.
To show injectivity suppose that pairs $(a:b)$ and $(a':b')$ define the same ideal class.
Then upon rescaling we can assume that $(a,b)=(a',b')$.
Now we have two surjective maps $R^2\to I=(a,b)$: one given by $(a,b)$ and another
by $(a',b')$, and our assertion follows from Lemma \ref{ideal-class-lem} below.
\ed

\begin{lem}\label{ideal-class-lem}
For every nonzero ideal $I\sub R$ 
the action of $\SL_2(R)$ on surjective maps $R^2\to I$ is transitive.
\end{lem}

\Pf . Since $I$ is a projective $R$-module, for every surjective
map $f:R^2\to I$ there exists an isomorphism
$$\a:R^2\rTo{\sim} I'\oplus I$$
such that $f$ is the composition of $\a$ with the projection to $I$.
Note that $\det(\a)$ induces an isomorphism of $R$ with $I'\ot_R I$, so
we obtain an isomorphism $I'\simeq I^{-1}$.
Thus, we can view $\a$ as an isomorphism
$R^2\to I^{-1}\oplus I$ such that $\det(\a)$ is the canonical isomorphism $R\to I^{-1}\ot I$.
If $g:R^2\to I$ is another surjective morphism and 
$\b:R^2\rTo{\sim} I^{-1}\oplus I$ is the corresponding isomorphism
then $\ga=\b^{-1}\circ\a$ is an element of $\SL_2(R)$ such that $g\circ\ga=f$.
\ed



\begin{rem}
In general the embedding \eqref{FM-map} is not a bijection as one can see already in
the case of a non-principally polarized abelian variety with $\End(A)=\Z$
(cf.\ \cite[Ex.\ 4.16]{Orlov-ab}). The Lagrangians not in the image of this map correspond
to categories of twisted sheaves equivalent to $D^b(A)$ (see \cite{P-sym}).
Note that the set $\LG_A(\Q)$ is a subset of vertices of
the spherical building associated with the group $\BU$, which is related
to the boundary of the Baily-Borel compactification of the Siegel domain $D_A$. 
It would be interesting to see whether other elements of this building have an interpretation in terms of $D^b(A)$. Also, one can expect some relation between the quasi-standard $t$-structures and the
$t$-structures associated with stabilities coming from points of $D_A$ or $D_A\times\C$.
In the case of K3-surfaces similar questions are studied in \cite{Ma} and \cite{Hart}.
\end{rem}




\begin{thebibliography}{99}
\bibitem{AP} D.\, Abramovich, A.\, Polishchuk, {\it Sheaves of t-structures and valuative
criteria for stable complexes}, Journal f\"ur die reine und angewandte Mathematik 590 (2006), 89--130. 
\bibitem{BMT} A.\, Bayer, E.\, Macr\`i, Y.\, Toda,
{\it Bridgeland Stability conditions on threefolds I: Bogomolov-Gieseker type inequalities},
arXiv:1103.5010.
\bibitem{BBD} A.~Beilinson, J.~Bernstein, P.~Deligne,
{\it Faisceaux pervers,} Ast\'erisque, 100, Soc. Math. France, Paris, 1982.
\bibitem{Block} J.\, Block, {\it Duality and equivalence of module categories in noncommutative geometry II: Mukai duality for holomorphic noncommutative tori}, arXiv:0604296.
\bibitem{Borel} A., Borel, {\it Density properties of certain subgroups of 
semisimple groups without compact components}, Ann. Math. 72 (1960), 179--188.
\bibitem{Borel-LAG} A.\, Borel, Linear algebraic groups, 2nd edition, Springer-Verlag, New York, 1991.
\bibitem{BoHC} A.\, Borel, Harish-Chandra, {\it Arithmetic subgroups of
algebraic groups}, Ann. Math. 75 (1962), 485--535.
\bibitem{Bridge-stab} 	  T.\, Bridgeland, 
	 {\it Stability conditions on triangulated categories}, Ann. Math. 166  (2007),  no. 2, 317--345.
\bibitem{Bridge-K3} T.\, Bridgeland, {\it Stability conditions on $K3$ surfaces},
Duke Math. J.  141  (2008),  no. 2, 241--291.
\bibitem{Douglas} M.\, R.\, Douglas, {\em  Dirichlet branes, homological
mirror symmetry, and stability},  Proceedings of the International
  Congress of Mathematicians, Vol. III (Beijing, 2002), 395--408, Higher
Ed. Press, Beijing, 2002. 
\bibitem{Giraud} J.\, Giraud, {\it Remarque sur une formula de Shimura-Taniyama},
Invent.~Math.~5 (1968), 231--236.
\bibitem{Go} R.\, Godement, {\it Domaines fondamentaux des groupes arithm\'etiques},
S\'eminaire Bourbaki, Vol. 8 (1962/63), Exp. 257, 201--225.
\bibitem{GLO} V.\, Golyshev, V.\, Lunts, D.\, Orlov, {\it Mirror symmetry for
abelian varieties}, J. Algebraic Geom. 10 (2001), 433--496.
\bibitem{Fulton} W.\, Fulton, Intersection theory, Springer-Verlag, 1984.
\bibitem{Hart} H.\, Hartmann, {\it Cusps of the K\"abler moduli space and stability conditions
on K3 surfaces}, arXiv:1012.3121.
\bibitem{HN} D.\, Huybrechts, M.\, Nieper-Wisskirchen,
{\it Remarks on derived equivalences of Ricci-flat manifolds}, arXiv:0801.4747.
\bibitem{Ma} S.\, Ma, {\it Fourier-Mukai Partners of a K3 Surface and the Cusps of its Kahler 
Moduli}, Internat. J. Math., 20 (2009), 727--750.
\bibitem{MMS} E.\, Macr\`i, S.\, Mehrorta, P.\, Stellari, {\it Inducing stability conditions}, preprint
math.AG/0705.3752.
\bibitem{Mukai-Fourier} S.\, Mukai, {\it Duality between $D(X)$ and $D(\hat{X})$ with its
application to Picard sheaves}.
Nagoya Math. J. 81 (1981), 153--175.
\bibitem{Mukai-bun} S.\, Mukai, {\it Semi-homogeneous vector bundles on an abelian variety},
J.~Math.~Kyoto Univ.\ 18 (1978), 239--272.
\bibitem{Mukai-spin} S.\, Mukai, 
{\it Abelian variety and spin representation} (in Japanese), in
{\it Hodge theory and algebraic geometry (Sapporo, 1994)}, English translation: Univ. of Warwick preprint, 1998.
\bibitem{Mum-ab} D.\, Mumford, Abelian varieties. Oxford Univ. Press, 1974.
\bibitem{Orlov-ab} D.\, Orlov, {\it 
Derived categories of coherent sheaves on abelian varieties and equivalences between them},
Izv. Math. 66 (2002), no. 3, 569--594.
\bibitem{P-thesis} A.\, Polishchuk,
{\it Biextensions, Weil representation on
         derived categories, and theta-functions}, Ph.D. Thesis,          Harvard University, 1996.
\bibitem{P-sym} A.\, Polishchuk,        {\it Symplectic biextensions and a generalization of the
         Fourier-Mukai transform},
         Math. Research Letters 3 (1996), 813--828.
\bibitem{P-maslov} A.\, Polishchuk,          {\it Analogue of Weil representation for abelian schemes},
         Journal f\"ur die reine und angewandte Mathematik 543 (2002), 1--37.
\bibitem{P-ab-var} A.\, Polishchuk,   Abelian Varieties, Theta Functions and the Fourier Transform,
          Cambridge University Press, Cambridge, 2003.
\bibitem{P-nc-tori} A.\, Polishchuk, {\it Classification of holomorphic vector bundles on 
          noncommutative two-tori},   
          Documenta Mathematica 9 (2004), 163--181.
\bibitem{P-const} A.\, Polishchuk,
{\it Constant families of t-structures on derived categories of coherent sheaves},
    Moscow Math. J. 7 (2007), 109--134.
\bibitem{P-ker-alg} A.\, Polishchuk,     {\it Kernel algebras and generalized Fourier-Mukai transforms}, J. Noncomm. Geom. 5 (2011), 153--251.
\bibitem{P-LIF} A.\, Polishchuk,  
{\it Lagrangian-invariant sheaves and functors for abelian varieties}, arXiv:1109.0527.
\bibitem{PS} A.\, Polishchuk, A.\, Schwarz,
{\it Categories of holomorphic vector bundles on noncommutative
          two-tori}, Comm. Math. Phys. 236 (2003), 135--159.  
\bibitem{Seidel} P.\, Seidel, {\it Graded Lagrangian submanifolds}, 
Bull. Soc. Math. France 128 (2000), no. 1, 103--149.
\bibitem{Weil}  A.\, Weil,   {\it   Sur   certains   groupes
d'op\'erateurs unitaires}, Acta Math. 111 (1964), 143--211.
\end{thebibliography}
\end{document}